\documentclass[12pt]{amsart}
\usepackage{amscd}
\usepackage{mathrsfs}
\usepackage{a4wide}
\usepackage{color}
\usepackage{amssymb,amsmath,amsthm,amsfonts,latexsym}
\usepackage{bbm} 
\usepackage{tikz}

\usepackage{amsmath,amsfonts,amssymb,amsthm} 
\usepackage{hyperref}

\begin{document}

\newtheorem{theorem}{Theorem}[section]
\newtheorem{lemma}[theorem]{Lemma}
\newtheorem{corollary}[theorem]{Corollary}
\newtheorem{claim}[theorem]{Claim}
\newtheorem{fact}[theorem]{Fact}
\newtheorem{affirmation}[theorem]{Affirmation}
\newtheorem{proposition}[theorem]{Proposition}
\newtheorem{example}[theorem]{Example}
\newtheorem{question}[theorem]{Question}
\newtheorem{definition}[theorem]{Definition}
\newtheorem{remark}[theorem]{Remark}
\newtheorem{notation}[theorem]{Notation}
\numberwithin{equation}{section}

\def\eqdef{\stackrel{\rm def}{=}}
\def\Proof{{\noindent {\em Proof:\ }}}
\def\fp{\hfill $\Box$}

\def\alex{\mbox{\sf AT}}
\def\baire{{\nat}^{\nat}}
\def\bairez{{\bf Z}^{\omega}}
\def\base{{\mathcal B} }
\def\binary{2^{< \omega}}
\def\bw{\sf{BW}}
\def\cantor{2^{\nat}}
\def\cantorx{2^X}
\def\ca{{\mathcal A}}
\def\cc{{\mathcal C}}
\def\ci{{\,mathcal I}}
\def\cj{{\mathcal J}}
\def\ck{{\mathcal K}}
\def\cs{{\mathcal S}}
\def\cR{{\mathcal R}}
\def\cL{{\mathcal{L}}}
\def\coA{{\mathbb{A}}}
\def\coB{{\mathbb{B}}}
\def\coL{{\mathbb{L}}}
\def\coK{{\mathbb{K}}}
\def\cl#1{\overline{#1}}
\def\cofin{\mbox{\sf CoFIN}}
\def\conv{\mbox{\sf conv}}
\def\ds{{\frak D}}
\def\esptop{(X,\top)}
\def\exh{\mbox{\sf Exh}}
\def\fsig{F_{\sigma}}
\def\fsd{F_{\sigma\delta}}
\def\fin{\mbox{\sf Fin}}
\def\filter{{\mathcal F}}
\def\fbw{\sf{FinBW}}
\def\afilter{\vec{\filter}}
\def\gfilter{\mathcal G}
\def\ged{G_\delta}
\def\hilbert{[0,1]^{\nat}}
\def\hm {O\check{c}(\mbox{$\tau$-closed, $\tau$})}
\def\ideal{{\mathcal I}}
\def\inte#1{\buildrel {\;\circ} \over #1}
\def\intseq{{\bf Z}^{<\omega}}
\def\infiN{\nat^{[\infty]}}
\def\ikl{[K,L]}
\def\iko{[K,O]}
\def\leqtau{\leq_{\tau}}
\def\k#1{{\mathcal K}(#1)}
\def\lx{<_{\scriptscriptstyle X}}
\def\ly{<_{\scriptscriptstyle Y}}
\def\mon{\sf{Mon}}
\def\rl{\Re_{l}}
\def\nat{\mathbb{N}}
\def\N{\mathbb{N}}
\def\nwd{\mbox{\sf nwd}}
\def\nulo{\mbox{\sf null}}
\def\opcla{{\bf cl}}
\def\nin{\not\in}
\def\oc#1{O\check{c}({\fin},{#1})}
\def\oppoint{$\mbox{\bf op}$}
\def\oqpoint{$\mbox{\bf oq}$}
\def\ptres{\emptyset\times\mbox{FIN}}
\def\prodxa{\prod_\alpha \xa}
\def\power#1{{\mathcal P}(#1)}
\def\powerx{{\mathcal P}(X)}
\def\ppoint{$\mbox{p}^+$}
\def\pp{\mbox{$\mathbf  p^+$}}
\def\Pm{\mbox{$\mathbf  p^-$}}
\def\qpoint{$\mbox{q}^+$}
\def\qq{\mbox{$\mathbf{q}^+$}}
\def\Q{\mathbb{Q}}
\def\qed{\hfill \mbox{$\Box$}}
\def\R{\mathbb{R}}
\def\rpoint{$\mbox{rp}^+$}
\newcommand{\lrangle}[1]{\langle#1\rangle}
\def\sd{\scriptstyle{\triangle}}
\def\seq{{\omega}^{<\omega}}
\def\seqN{\nat^{<\omega}}
\def\sfan{S(\omega)}
\def\somega{S_\omega}
\def\soi{S(\ideal)}
\def\su{\subseteq}
\def\sumi{\mbox{\sf Sum}}
\def\swc{\overline{S_{\Omega}}}
\def\Sq{\mbox{$\mathbf{p}^-$}}
\def\sw{S_{\Omega}}
\def\sx{\scriptscriptstyle X}
\def\sy{\scriptscriptstyle Y}
\def\tauf{\tau_{\scriptscriptstyle \filter}}
\def\tauaf{\tau_{\scriptscriptstyle \afilter}}
\def\tauc{\tau_{\scriptscriptstyle \cc}}
\def\tauv{\tau_{V}}
\def\ti{T(\ideal)}
\def\topy{\top_{\sy}}
\def\topo{topolog\II a }
\def\top{\tau}
\def\ult{\mathcal U}
\def\vkl{{\mathcal V}_{\scriptscriptstyle{K,L}}}
\def\veba#1#2{{\mathcal V}_{\scriptscriptstyle{#1,#2}}}
\def\iveba#1#2{[#1,#2]}
\def\ve{{\mathcal V}}
\def\ufilter{\mathcal U}
\def\w1{\omega_1}
\def\wpp{\mbox{$\mathbf{wp}^+$}}
\def\ws{\mbox{$\mathbf{ws}$}}
\def\xa{X_\alpha}
\def\Z{\mathbb{Z}}

\def\sequen{\nat^{<\bf \scriptscriptstyle \omega}}
\def\casig{=^*}

\def\pie{\Pi_1^1}
\def\pieb{\boldsymbol{ \Pi}_1^1}
\def\pie03{{\mathbf{\Pi}}_3^0}
\def\sigb{{\mathbf \Sigma}_1^1}
\def\delb{\boldsymbol{ \Delta}_1^1}
\def\si12{\Sigma^{1}_{2}}
\def\del{\Delta_1^1}
\def\zig{\Sigma_1^1}
\def\pai2{\mathbf{\Pi}_2^0}
\def\sig2{{\mathbf{\Sigma}}^{0}_{2}}

\title{Ideals on countable sets: a survey with questions}
\author{Carlos Uzc\'ategui Aylwin}

\address{Escuela de Matem\'aticas, Facultad de Ciencias, Universidad Industrial de
	Santander, Ciudad Universitaria, Carrera 27 Calle 9, Bucaramanga,
	Santander, A.A. 678, COLOMBIA. Centro Interdisciplinario de L\'ogica y \'Algebra, Facultad de Ciencias, Universidad de Los Andes, M\'erida, VENEZUELA.}
\email{cuzcatea@saber.uis.edu.co}

\date{\today}

\thanks{The author thank La Vicerrector\'ia de Investigaci\'on y Extensi\'on de la Universidad Industrial de Santander for the financial support for this work,  which is part  of the VIE project  \#2422}



\begin{abstract}
An ideal on a set $X$ is a collection of subsets of $X$ closed under the operations of taking  finite unions and  subsets of its elements.  Ideals are  a very useful notion in topology and set theory and have been studied for a long time. 
We present a survey of results  about ideals on countable sets and include many open questions. 
\end{abstract}

\maketitle

\section{Introduction}

An ideal on a set $X$ is a collection of subsets of $X$ closed under the operations of taking  finite unions and  subsets of its elements.  Ideals are  a very useful notion in topology and set theory and have been studied for a long time. 
We present a survey of results  about ideals on countable sets and include many open questions.

We have tried to include aspects that were not covered in the survey written by  M. Hru{\v{s}}{\'a}k \cite{Hrusak2011}.  
We start by presenting two common forms to  define ideals: based on submeasures or on collections of nowhere dense sets. A basic tool in the study of ideals  are some orders to compare them: Kat\v etov, Rudin-Keisler and Tukey order. We focus mostly on the Kat\v etov order. 
The reader can consult \cite{Hrusak2011,solecki2015,solectodor2011,SolecTodor2004} for results on the Tukey order.  One important ingredient of our presentation is that we deal mainly  with definable ideals: Borel, analytic or co-analytic ideals. Another  crucial aspect is the role played by combinatorial properties of ideals,  a theme that has been  very much studied and provides a common ground  for the whole topic. Most of the work on ideals has been concentrated on tall ideals, nevertheless we include a section on Fr\'echet ideals (i.e.,  locally non tall ideals).  Since the properties about ideals we are dealing with are, in one way or another, based on selection principles, we end the paper with a discussion of Borel selection principles for ideals, that is,  the selection function is required to be Borel measurable.

We do not pretend to give  a complete revision of this topic; in fact,  the literature is vast and we have covered a small portion of it. Our purpose was to present some of the diverse ideas that have being used for studying ideals on countable sets and collect some open questions  which were scattered in the literature. 

\section{Terminology}
\label{prelimi}

An ideal  $\ideal$  on a set $X$ is a collection of subsets of $X$ such that:
\begin{itemize} 
\item[(i)] $\emptyset \in \ideal$ and $X\nin \ideal$. 
\item[(ii)] If $A, B\in \ideal$, then $A\cup B\in \ideal$. 
\item[(iii)] If $A\su B$ and $B\in \ideal$, then $A\in \ideal$. 
\end{itemize}
Given an ideal $\ideal$ on $X$, the {\em dual filter} of $\ideal$, denoted $\ideal^*$,  is the collection of all sets $X\setminus A$ with $A\in \ideal$.  We denote by  $\ideal^+$ the collection of all subsets of  $X$  which do not belong to $\ideal$.
Two ideals $\ideal$ and $\cj$ on $X$ and $Y$ respectively are {\em isomorphic} if there is a bijection $f:X\to Y$ such that $E\in \ideal$ if, and only if, $f[E]\in \cj$.  Suppose $X$ and $Y$ are disjoint, then the {\em free sum} of $\ideal$ and $\cj$, denoted by $\ideal\oplus\cj$ is defined on $X\cup Y$ as follows: $A\in \ideal\oplus \cj$ if $A\cap X\in \ideal$ and $A\cap Y\in \cj$. 

We denote by $2^{<\omega}$  (respectively, $\seqN$) the collection of all finite binary sequences (respectively, finite sequences of natural numbers). If $x\in \cantor$, then $x\restriction n$ is the sequence $\langle x(0),\cdots, x(n-1)\rangle$ for $n\in \nat$.

Now we recall some combinatorial properties of ideals.  We put $A\su^*B$ if $A\setminus B$ is finite.
An ideal $\ideal$ is a $P$-ideal, if for any family $E_n\in \ideal$ there is $E\in \ideal$ such that $E_n\su^* E$ for all $n$.  This is one of the most studied class of ideals. 

\begin{enumerate}
\item[({$\mathbf  p^+$})] $\ideal$ is \pp, if for every decreasing sequence $(A_n)_n$ of sets in $\ideal^+$, there is $A\in \ideal^+$ such that $A\su^* A_n$ for all
$n\in\nat$. Following  \cite{HMTU2017}, we say that $\ideal$   is $\Pm$, if for every decreasing sequence $(A_n)_n$ of sets in $\ideal^ +$ such that $A_n\setminus A_{n+1}\in \ideal$, there is $B\in \ideal^+$ such that $B\su^* A_n$ for all $n$. 

The following notion was suggested by some results  in \cite{Farah2003,Filipowetal2008}.  Let us call a scheme a collection $\{A_s: \; s\in\binary\}$ such that $A_s=A_{s\widehat{\;\;}0}\;\cup\; A_{s\widehat{\;\;}1}$ and $A_{s\widehat{\;\;}0}\cap A_{s\widehat{\;\;}1}=\emptyset$ for all
$s\in\binary$.  An ideal  is \wpp,  if for every scheme $\{A_s:\; s\in\binary\}$ with $A_\emptyset\in \ideal^+$, there is $B\in\ideal^+$ and $\alpha\in\cantor$ such that $B\su^* A_{\alpha\restriction n}$
for all $n$. 

\item[(\qq)] $\ideal$ is \qq, if for every $A\in \ideal^+$ and every
partition $(F_n)_n$ of $A$ into finite sets, there is $S\in\ideal^+$
such that $S\su A$ and $S\cap F_n$  has at most one element for
each $n$. Such sets $S$ are called (partial) {\em selectors} for the partition. If we allow partitions with pieces in $\ideal$, we say that the ideal is {\em weakly selective} \ws\  \cite{HMTU2017} (also called weakly Ramsey in 
\cite{SametTsaban2009}). Another natural variation is as follows: For every partition $(F_n)_n$ of a set $A\in\ideal^+$ with each piece $F_n$ in $\ideal$,  there is $S\in\ideal^+$ such that $S\su A$ and  $S\cap F_n$ is finite for all $n$. It is known that the last property  is  equivalent to \Pm\ (see  Theorem \ref{ppoint}).
\end{enumerate}

All spaces are assumed  to be regular and $T_1$.   A collection $\base$ of non empty open subsets  of $X$ is a {\em $\pi$-base}, if every non empty open set contains a set belonging to $\base$.  A point $x$ of a topological space $X$ is called a {\em Fr\'echet
point}, if for every $A$ with $x\in \cl A$  there is a sequence $(x_n)_n$ in $A$
converging to $x$.   It is well known that filters (or dually, ideals) are viewed as spaces with only one non isolated point. We recall this basic construction. Suppose $Z=\nat\cup \{\infty\}$ is a space such that  $\infty $ is the only accumulation point.  Then  $\mathcal{F}_\infty=\{A\su \nat:\; \infty \in \mathrm{int}_Z(A\cup \{\infty\} )\}$ is  the neighborhood filter of $\infty$. 
Conversely, given an ideal $\mathcal{I}$ over $\nat$,  we define  a topology on $\nat\cup \{\infty\}$  by declaring that  each $n\in \nat$ is isolated and $\mathcal{I}^*$ is the neighborhood filter of $\infty$. We denote this space by $Z(\mathcal{I})$.    It is clear that the combinatorial properties of $\ideal$ and $Z(\ideal)$ are the same.

For $n\in\nat$, we denote by $X^{[n]}$ the collection of  $n$-elements subsets of $X$ and $X^{[\infty]}$ the collection of infinite subsets of $X$. 
The classical  Ramsey theorem asserts that for every coloring
$c:\nat^{[2]}\to \{0,1\}$, there is an infinite subset $X$ of $\nat$ such that $X$
is \emph{$c$-homogeneous}, that is,  $c$ is constant in $X^{[2]}$. 
An ideal $\ideal$  is  \emph{Ramsey at} $\nat$,    when for any coloring $c:\nat^{[2]}\to \{0,1\}$ there is a $c$-homogeneous set which is $\ideal$-positive,  we denoted it by $\nat\to(\ideal^+)^2_2$. If it is the case that for any coloring $c$ and any $X\in\ideal^+$ there is a  $c$-homogeneous set $Y\in \ideal^+$ contained in $X$, we shall write   $\ideal^+\to(\ideal^+)^2_2$ and called such ideal a {\em Ramsey} ideal.
A collection $\ca$ of subsets of a set $X$ is {\em tall}, if for every infinite set $A\su X$, there is an infinite set $B\su A$ with $B\in \ca$.  Ramsey's theorem says that the collection of $c$-homogeneous sets is a tall family for any coloring $c$.

A general reference for all descriptive set theoretic notions used in this paper is \cite{Kechris94}. A set is $\fsig$ (also denoted  $\sig2$)  if it is equal to  the union of a countable collection of closed sets. Dually, a set is $\ged$ (also denoted $\pai2$) if it is the intersection of a countable collection of open sets.  The Borel hierarchy is the collection of classess $\mathbf{\Sigma}^0_\alpha$ and $\mathbf{\Pi}^0_\alpha$ for $\alpha$ a countable ordinal. For instance,   $\pie03$ (which is also denoted by $F_{\sigma\delta}$) are the sets of the form $\bigcap_n  F_n$ where each $F_n$ is an $\fsig$. 
 A subset $A$ of a Polish space is called {\em analytic}, if it is a
continuous image of a Polish space. Equivalently, if there is a
continuous function $f:\baire\rightarrow X$ with range $A$, where
$\baire$ is the space of irrationals.    Every Borel subset of a Polish space is analytic. A subset of a Polish space  is {\em co-analytic} if its complement is analytic. The class of analytic (resp., co-analytic) sets is denoted by $\sigb$ (resp. $\pieb$).

\section{Some examples}
\label{ejem}

In this section we present some examples of ideals. The interested reader can consult \cite{Hrusak2011,HMTU2017,meza-tesis}  where he can find many more interesting examples. 

The simplest ideal is $\fin$, the collection of all finite subsets of $\nat$. There are two natural ideals quite related to $\fin$ which are defined on $\nat\times \nat$. 

\begin{equation}
\label{vacioporfin}
\{\emptyset\}\times \fin =  \{A\su\nat\times\nat:\; \{m\in\nat:\; (n,m)\in A\}\in\fin\;\mbox{for all $m\in\nat$}\}.
\end{equation}

\begin{equation}
\label{finporvacio}
\fin\times\{\emptyset\} =  \{A\su\nat\times\nat:\; \exists n\in\nat\; A\su \{0,1,\cdots, n\}\times \nat\;\}. 
\end{equation}

In general, let $\ideal$ and $\cj$ be ideals on $X$ and $Y$ respectively; its {\em Fubini product} $\ideal\times\cj$ is an ideal on $X\times Y$ defined as follows: for $A\su X\times Y$, we let $A_x=\{y\in Y\; (x,y)\in A\}$. 
\[
\ideal\times\cj=\{A\su X\times Y:\; \{x\in X:\; A_x\nin \cj\}\in \ideal\}.
\]
By an abuse of notation,  the ideals $\{\emptyset\}\times \fin$  and $\fin\times \{\emptyset\}$ are usually denoted $\emptyset\times\fin$ and $\fin\times \emptyset$, respectively.  An ideal $\ideal$ on $X$ is {\em countably generated} if there is a countable collection $\{A_n:\; n\in\nat\}$ of subsets of $X$ such that $E\in \ideal$ if, and only if, there is $n$ such that $E\su A_0\cup\cdots\cup A_n$.  The only countably generated ideals containing  all finite sets are $\fin$ and $\fin\times\{\emptyset\}$ (see Proposition 1.2.8. in \cite{Farah2000}).

Two very important ideals on $\Q$ are the ideal of nowhere dense subsets of $\Q$ (with its usual metric topology), denoted $\nwd(\Q)$, and the ideal of null sets defined as follows:
\[
\nulo(\Q)  =  \{A\su \Q\cap [0,1]:\; \cl{A}\;\mbox{has Lebesgue measure zero}\}.
\]
In general, if $X$ is a topological space, then $\nwd(X)$ denotes the ideal of nowhere dense subsets of $X$.  Another very natural ideal associated to a space is defined as follows:
For every  point $x\in X$, let 
\begin{equation}
\label{acumIdeal}
 \cj_x=\{A\su  X:\; x\nin \cl{A\setminus\{x\}}\}.
\end{equation}
In fact, every ideal on $X$ is of the form $\cj_x$ for some topology on $X$.

Two ideals on $\nat$ that have a very natural connection with number theory and real analysis are the following: 
\[
\ideal_{1/n}=\{A\su\nat:\;\sum_{n\in A}\frac{1}{n+1}<\infty\}
\]
and 
\[
A\in\ideal_d\;\;\Leftrightarrow \;\; \limsup_n \frac{|A\cap\{0,1,\cdots, n-1\}|}{n}=0.
\]
The ideal  $\ideal_d$ consists of the {\em asintotic density zero sets}.

Let $CL(\cantor)$ denote the collection of clopen subsests of $\cantor$ and $\lambda$ the product measure on $\cantor$.  Notice that $CL(\cantor)$ is countable. Let
\[
\Omega=\{A\in CL(\cantor):\; \lambda(A)=1/2\}.
\]
Solecki's ideal  $\cs$ is the ideal on $\Omega$ generated by the following sets: 
\[
E_x=\{A\in \Omega:\; x\in A\},
\]
where  $x\in \cantor$.   Solecki introduced $\cs$ to characterize the ideals satisfying Fatou's lemma \cite{Solecki2000}.

The following ideal is called the \emph{eventually different} ideal:
$$
\mathcal{ED}=\{A\subseteq\nat\times\nat: (\exists n)(\forall m\geq n)(|{\{(m,k):\; k\in\nat\}\cap A}|\leq n)\}.
$$
The following restriction of $\mathcal{ED}$ also plays an important role in the study of combinatorial properties of ideals: 
$$
\mathcal{ED}_{\fin}=\mathcal{ED}\restriction \Delta,
$$
where $\Delta=\{(n,m)\in\nat\times\nat:\;m\leq n\}$. Note that $\mathcal{ED}_{\fin}$ is (up to isomorphism)
the unique ideal generated by the selectors of some partition of
$\nat$ into finite sets $\{I_n:n\in\nat\}$ such that
$\limsup_n |{I_n}|=\infty$.  As we will see later,    $\mathcal{ED}_{\fin}$ is critical for the
\qpoint-property.

Let $\conv$ be the ideal generated by the range of all convergent sequences of rationals numbers, where the convergence is in  $\R$. In other words, $\conv$ is the collection of all subsets of $\Q$ such that  the Cantor-Bendixon derivative of its closure in $\R$ is finite. 

\bigskip

Now we present a family of ideals defined by homogeneous sets for colorings.  Let $c:\nat^{[2]}\to \{0,1\}$ be a coloring. Recall that a set $H\su\nat$ is $c$-homogeneous if $c$ is constant in $H^{[2]}$. The collection $hom(c)$ of all $c$-homogeneous sets is closed in $\cantor$.  Let $\ideal_{hom(c)}$ be the ideal generated by the $c$-homogeneous sets.

The \emph{infinite random graph} on $\nat$, also known as  the \emph{Rado graph}  or the  \emph{Erd\"os-R\'enyi graph} (see, e.g. ,\cite{Cameron1997}) can be concisely described as
follows. Recall that a family $\{X_n:\; n\in\nat\}$ of infinite subsets of $\nat$ is \emph{independent}, if given two disjoint finite subsets $F, E$ of $\nat$ the set $(\bigcap_{n\in F}X_n ) \setminus (\bigcup_{n\in E}X_n)$ is infinite.
Let $\{X_n:\, n\in\nat\}$ be an independent family of subsets of $\nat$ such that 
$n\in X_m$ if, and only if, $m\in X_n$, for all $m,n\in\nat$. The random graph is then $(\nat, E)$, where
$$
E=\{\{n,m\}:m\in X_n\}.
$$
The random graph is \emph{universal} in the following sense.  Given a graph
$\langle{\nat,G}\rangle$, there is a subset $X\subseteq\nat$ such
that $\lrangle{\nat,G}\cong \lrangle {X, E\restriction X}$.  
The \emph{random graph ideal} $\mathcal{R}$ is the ideal on $\nat$ generated by cliques and free sets  of the random graph or, equivalently,  the homogeneous sets with respect to the \emph{random coloring} $c:[\nat]^2\to 2 $ defined by  $c (\{m,n\})=1$ if, and only if, $\{m,n\}\in E$.

\section{Complexity of ideals}
\label{complexity}

 We say
that a collection $\ca$ of subsets of  a countable set $X$ is {\em analytic} (resp. Borel),
if $\ca$ is analytic  (resp. Borel) as a subset of the cantor cubet $2^X$
(identifying subsets of $X$ with characteristic functions) \cite{Kechris94}. The set $\infiN$ of infinite subsets of $\nat$ will be always considered with the subspace topology of $\cantor$.  We say that an ideal is {\em analytic}, if it is an analytic as a subset of $2^X$.  Since the collection of finite subsets of $\nat$ is a dense set in $\cantor$, then there are no ideals containing $\fin$ which are closed as subsets of $\cantor$.  On the other hand, if $\ideal$ is a $G_\delta$ ideal with $\fin \su \ideal$, then $
\ideal^*=\{\nat\setminus A:\; A\in \ideal\}$ is also dense $\ged$ (as the map $A\mapsto \nat\setminus A$ is an homeomorphism of $\cantor$ into itself). Therefore by the Baire category theorem  $\ideal \cap \ideal^*\neq \emptyset$, which says that $\ideal =\power\nat$.  So,  the simplest Borel   ideals have complexity $\fsig$. They have been quite investigated as we will see.  
Every analytic ideal is generated by a $\ged$ set, i.e., there is a $\ged$  $G\su \ideal$ such that every $A\in \ideal$ is a subset of a finite union of elements of $G$  \cite{Zafr90} (see also \cite[Theorem 8.1]{SolecTodor2004}). 

Most of the theory of definable ideals has been concentrated on analytic ideals. There are a few results about co-analytic ideals. 
The following theorem provides a very general representation of  analytic ideals  on spaces of  continuous functions.  It is an instance of the ideal defined by \eqref{acumIdeal}. 

\begin{theorem}
\label{repre1}
(Todor\v{c}evi\'c \cite[Lemma 6.53]{Todor2010})
Let $\ideal$ be an ideal over $\nat$. The following are equivalent.
\begin{enumerate}
\item[(i)] $\ideal$ is analytic.

\item[(ii)] There  are  continuous functions $f, f_n:\baire\to \R$ for $n\in \nat$ such that  $f$ is an accumulation point of $\{f_n:\; n\in\nat\}$ respect to product topology on $C_p(\baire)$ and 
\[
E\in \ideal \;\;\mbox{if, and only if, }  f\nin \overline{\{f_n:\; n\in E\}},
\]
where the closure is taken in $C_p(\baire)$. 
\end{enumerate}

\end{theorem}

There are some well known co-analytic ideals. Let $WO(\Q)$ be the ideal of well ordered subsets of $\Q$. This is a typical complete co-analytic  ideal. 
Let $\ideal_{wf}$ be the ideal on $\seqN$ generated by the well founded
trees on $\N$, i.e., $A\su \N^{<\omega}$ belongs to $\ideal_{wf}$,  if there is a
wellfounded tree $T$ such that $A\su T$. This is
equivalent to say that the tree generated by $A$ is well founded.
Then $\ideal_{wf}$ is also a complete co-analytic ideal.  In sections  \ref{toprepre} and \ref{frechet} we shall present another examples of co-analytic ideals (see also   \cite{Farkasetal2018,meza-tesis}).

We do not know of any general theorem,  as Theorem \ref{repre1},  for co-analytic ideals.  So we state this question as follows.

\begin{question}
Is there a general representation theorem for co-analytic ideals?
\end{question}

\section{Ideals based on submeasures}

A natural and very impotant method for defining ideals  is based on  measures or, more generally, submeasures. In this section we present some of these ideas. 

A function  $\varphi:\mathcal{P}(\nat)\to[0,\infty]$ is a {\em lower semicontinuous submeasure (lscsm)} if $\varphi(\emptyset)=0$,
$\varphi(A)\leq\varphi(A\cup B)\leq\varphi(A)+\varphi(B)$ and $\varphi(A)=\lim_{n\to\infty}\varphi(A\cap\{0,1,\cdots, n\})$.

There are several ideals associated to a lscsm: 
\[
\begin{array}{lcl}
\fin(\varphi)& =& \{A\subseteq\nat:\varphi(A)<\infty\}.
\\
\exh(\varphi)& = &\{A\subseteq\nat: \lim_{n\to\infty}\varphi(A\setminus\{0,1,\dots, n\})=0\}.
\\
\sumi(\varphi) &=&\{A\subseteq \mathbb{N}:\sum_{n\in A}\varphi\{n\}<\infty\}.
\end{array}
\]
They satisfy the following relations:
$$
\sumi(\varphi)\subseteq \exh(\varphi)\subseteq \fin(\varphi).
$$
The collection of ideals that can be represented by one these three forms have been extensively investigated. The work of Farah \cite{Farah2000} and Solecki \cite{Solecki1999} are two of the most important early works for the study of the ideals associated to submeasures.  

To each divergent series  $f: \nat\to [0,+\infty)$   of possitive real numbers,  we associate a measure on $\nat$ by 
\[
\varphi_f(A)=\sum_{n\in A}f(n).
\]
An ideal $\ideal$ is \emph{summable}  \cite{Mazur91} if there is a divergent series $f$ as above such that $\ideal=\fin(\varphi_f)$. Notice that $\sumi(\varphi_f)=\fin(\varphi_f)$. The  usual notation for this ideal is $\ideal_f$.  A typical example is the following 
\[
\ideal_{1/n}=\{A\su\nat:\;\sum_{n\in A}\frac{1}{n+1}<\infty\}.
\]

Another very natural way of defining lscsm is as follows. Let $I_n$ be a partition of $\nat$ into finite sets. Let $\nu_n$ be a measure
on $I_n$ (i.e., there is a function $f_n:I_n\to [0,+\infty)$ such that $\nu_n(A)=\sum_{n\in A}f_n(n)$). Let 
\[
\varphi(A)=\sup_n v_n(A\cap I_n).
\]
Then $\varphi$ is a lscsm and $Exh(\varphi)$ is called a {\em density ideal} \cite{Farah2000}.  The prototype is the following 

\begin{example} Let $\varphi_d:\power\nat\to [0,+\infty]$ given by 
\[
\varphi_d(A)=\sup\left\{\frac{A\cap \{0,1,\cdots, n-1\}}{n}:\; n\in \nat\right\}.
\]
Then $\ideal_d=\exh(\varphi_d)$ is the ideal of asintotic density zero sets.  We have
\[
A\in\ideal_d\;\;\Leftrightarrow \;\; \limsup_n \frac{|A\cap\{0,1,\cdots, n-1\}|}{n}=0.
\]
\end{example}

The Cantor set $\{0,1\}^\nat$ is a group with the product operation where $\{0,1\}$ is the  group $\Z_2$; equivalently, viewing the elements of $\cantor$ as subsets of $\nat$, then the algebraic operation is the symmetric difference.  Then $\cantor$ is actually a  Polish group.  Every ideal on $\nat$ is a subgroup of $\cantor$. Since there are no $\ged$ ideals (containing $\fin$), then none of these subgroups are Polish. However, the following weaker notion has been used to study subgroups of Polish groups. We say that a subgroup $G$ of $\cantor$ is {\em Polishable}, if there is a Polish group topology on $G$ such that the Borel structure of this topology is the same as the Borel structure $G$ inherites from $\cantor$. 

The following representation  of analytic $P$-ideals is the most fundamental result about them.  It says that any $P$-ideal is in a sense similar to a density ideal. 

\begin{theorem}
\label{solecki}
(S. Solecki \cite{Solecki1999}) Let $\ideal$ be an analytic ideal on $\nat$. The following are equivalent: 
\begin{itemize}
\item[(i)]  $\ideal$ is a $P$-ideal.
\item[(ii)] There is a lscsm $\varphi$ such that $\ideal=\exh(\varphi)$. 
\item[(iii)] $\ideal$ is Polishable. 
\end{itemize}
In particular, every analytic $P$-ideal  is $F_{\sigma\delta}$.  Moreover, $\ideal$ is an $\fsig$ $P$-ideal, if,  and only if, there is a lscsm $\varphi$ such that $\ideal=\exh(\varphi)= \fin(\varphi)$.
\end{theorem}

\subsection{$\fsig$ and $F_{\sigma\delta}$ ideals.}
\label{fsig-fsigdelta}

As we said, from the complexity point of view, $\fsig$ ideals are the simplest ones. In this section we present some results about them. 

 A set $\ck \su\power\nat$ is {\em hereditary} if for every $A\in \ck$ and  $B\su A$ we have that $B\in \ck$.  A family $\ck$ of subsets of $\nat$ is said to be {\em closed under finite changes} if $A{\sd} F\in \ck$ for every $A\in \ck$ and  a finite set $F\su \nat$. Given an hereditary collection $\ck$, we denote by $\ideal_\ck$ the ideal generated by $\ck$. That is to say
\[
B\in \ideal_\ck\;\;\Leftrightarrow \;\; \exists n\in \nat, \exists A_1,\cdots, A_n \in \ck\; (B\su A_1\cup\cdots\cup A_n).
\]

\begin{theorem} (Mazur \cite{Mazur91})
Let $\ideal$ be an ideal on $\nat$. The following are equivalent.
\begin{itemize} 
\item[(i)] $\ideal$ is $\fsig$.
\item[(ii)] there is a hereditary  closed collection $\ck$ of subsets of $\nat$ such that $\ideal=\ideal_\ck$.

\item[(iii)] there is a lscsm $\varphi$ such that $\ideal=\fin(\varphi)$.

\end{itemize}

\end{theorem}

An important example of $\fsig$ ideals are  $\ideal_{hom(c)}$, the ideal generated by the family $hom(c)$ of homogeneous sets respect to a coloring $c$ (see \S \ref{ejem}).   Notice that $hom(c)$ is a closed hereditary collection of subsets of $\nat$. 

The following is  part of the folklore  (for a proof see e.g., \cite[Lemma 3.3]{HMTU2017}).

\begin{theorem}
\label{fsigpmas}
Every $\fsig$ ideal is \ppoint. 
\end{theorem}

An important question involving  $\fsig$ ideals is the following:

\begin{question} 
(M. Hru\v s\'ak, \cite{Hrusak2017}) Does every tall Borel ideal contain a tall $\fsig$ ideal?
\end{question}

The previous question can be understood as asking whether  an analog of the classical perfect set theorem holds for the collection of  tall Borel ideals. However, the analogy is not complete, since there exists a   $\mathbf{\Pi}^1_2$  ideal which does not contain any $\fsig$ tall ideal (see \cite[Theorem 4.24]{Grebikuzca2018}).

We have seen in Theorem \ref{solecki} that every analytic $P$-ideal is $F_{\sigma\delta}$. One could naturally ask whether such ideals are a countable intersection of $\fsig$ ideals. Since this is not true in general,  Farah  \cite{farah2004luzin} introduced a weaker property (we follow the presentation given in \cite{hrusakmeza2011}).  They called an ideal $\ideal$ {\em Farah} if there is a countable collection  $\ck_n$ of closed hereditary families of subsets of $\nat$ such that 
\[
\ideal=\{A\su\nat:\, (\forall n\in \nat)(\exists m\in \nat)(A\setminus\{0,1,\cdots, m\}\in \ck_n)\}.
\]
It is clear that every Farah ideal is $F_{\sigma\delta}$. 
In \cite{farah2004luzin} it is shown that $\nwd(\Q)$,  $\nulo(\Q)$ and every analytic $P$-ideal are  Farah. However, there is no an $\fsig$ ideal $\cj$ such that $\nwd(\Q)\su \cj$. 

\begin{theorem}
(M. Hru\v s\'ak and D. Meza-Alc\'antara, \cite{hrusakmeza2011})  Let $\ideal $ be an ideal on $\nat$. The following are equivalent:
\begin{itemize}
\item[(i)] $\ideal$ is Farah. 
\item[(ii)] There is a sequence $\{F_n:\; n\in \nat\}$ of hereditary $\fsig$ sets closed under finite changes such that  $\ideal=\bigcap_n F_n$.
\item[(iii)] There is a sequence $\{F_n:\; n\in \nat\}$ of $\fsig$ sets closed under finite changes such that  $\ideal=\bigcap_n F_n$.
\end{itemize}
\end{theorem}

The previous result suggests a weaking of the notion of a Farah ideal. An ideal $\ideal$ is called {\em weakly Farah} \cite{hrusakmeza2011} if there is a sequence $\{F_n:\; n\in \nat\}$ of hereditary $\fsig$ sets such that  $\ideal=\bigcap_n F_n$. 

\begin{question}
(i)  (Farah \cite{farah2004luzin}) Is every $F_{\sigma\delta}$ ideal a Farah ideal?

(ii) (M. Hru\v s\'ak and D. Meza-Alc\'antara, \cite{hrusakmeza2011}) Is every $F_{\sigma\delta}$ ideal weakly Farah? Is every weakly Farah ideal a Farah ideal?
\end{question}

\subsection{Summable ideals on Banach spaces.}
\label{summableBanach}

The notion of a summable ideal has been extended to ideals where the sum is calculated in a Banach space or, more generally, in a Polish abelian group \cite{Borodulinetal2015,Boban1999,Farah1999}.  In this section we present some results and questions about this approach.

Let $G$ be a Polish abelian group (with additive notation) or a Banach space. Let $h:\nat\to G$ be a sequence. We say that the series $\sum_n h(n)$ is {\em unconditional convergent} in $G$ if  the net $\{\sum_{n\in A}h(n):\; A\in \fin \}$ (where $\fin$ is ordered by $\subseteq$) is convergent in $G$. This is equivalent to requiere that $\sum_n h(\pi(n))$ is convergent in $G$ for every permutation $\pi$ of $\nat$. Let $h:\nat\to G$ be a sequence such that $\sum_n h(n)$ does not exist. The {\em generalized summable ideal} $\ideal^G_h$ associated to $G$ and $h$ is the following \cite{Borodulinetal2015}:
\[
A\in \ideal^G_h \;\;\Leftrightarrow\;\; \sum_{n\in A}h(n)\;\;\mbox{is unconditional convergent in $G$}.
\]
An ideal $\ideal$ is said to be $G$-representable, if there is $h$ such that $\ideal=\ideal^G_h$. Analogously, it is defined when an ideal is $\cc$-representable for $\cc$ a class of abelian Polish groups. 

We recall that a lscsm $\varphi$ is {\em non-pathological} \cite{Farah2000} if $\varphi(A)$ is equal to the supremum of all $\nu(A)$ for $\nu$  a measure such that $\nu\leq \varphi$. A $P$-ideal is {\em non-pathological}, if it is equal to $Exh(\varphi)$ for some non-pathological lscsm $\varphi$.

\begin{theorem}(Borodulin-Nadzieja, Farkas, Plebanek \cite{Borodulinetal2015})
\begin{itemize}
\item[(i)] An ideal is $\R$-representable if, and only if, it is summable. 

\item[(ii)] An ideal is Polish-representable if, and only if, it  is an analytic $P$-ideal. 

\item[(iii)] An analytic $P$-ideal is Banach-representable if, and only if, it is non-pathological. 

\item[(iv)] A tall $\fsig$ $P$-ideal is representable in $c_0$ if, and only if, it is summable.

\item[(v)] There is an $\fsig$ tall ideal representable in $l_1$ which is not summable.
\end{itemize}
\end{theorem}

\begin{question}\cite{Borodulinetal2015}
How to characterize analytic $P$-ideals which  are $c_0$-representable?
\end{question}

\begin{question}\cite{Borodulinetal2015}
How to characterize ideals which are $l_1$-representable? Are they necessarily  $\fsig$?
\end{question}

\section{Topological representations by nowhere dense sets}
\label{toprepre}

In this section we review some constructions of ideals motivated by the ideal of nowhere dense sets. We consider two  different ways of presenting $\nwd(\Q)$.  For the first one, we see $\Q$ as a dense subset of $\R$ and we have  the following representation:
$$
\nwd(\Q)  =  \{A\su \Q\cap [0,1]:\; \cl{A}\;\mbox{is nowhere dense in $\R$}\}.
$$
On the other hand,  if $\base$ is a base for $\Q$ (of non empty open sets), we  have
\[
A\in \nwd(\Q)\;\; \Leftrightarrow \;\; \forall V    \in \base, \exists W\in \base \;(W\su V\wedge W\cap A=\emptyset).
\]
We will address each of these approaches in this section.

\subsection{Ideals of nowhere dense sets.}
A natural question is to determine for a given ideal $\ideal$ on a set $X$ whether there is a topology $\tau$ on $X$ such that $\ideal=\nwd(\tau)$ . This question was studied in \cite{CieJas95}, but most of their results are for $X$ uncountable. For $X$ countable, in \cite{todoruzca} are shown some general negative results (i.e., ideals for which such topology does not exist).

 Since we are mostly interested in definable ideals, we will work with {\em analytic topologies}, i.e., topologies $\tau$ on $X$ such that $\tau$ is analytic as a subset of $2^ X$ (see \S\ref{complexity}). The study of analytic topologies was initiated in  \cite{todoruzca}  (see also \cite{CamargoUzca2018b,CamargoUzca2018,Shiva2017,Todoruzca2000,Todoruzca2014,Uzca03,UzVi05}). 

Let $\mathcal I$ be an ideal over $X$ containing all singletons. Then the
dual filter (together with $\emptyset$) is a $T_1$ (but not Hausdorff)
topology  such that its nowhere dense sets are exactly the sets in  $\mathcal I$.
The next natural question is to requiere that the topology is $T_2$. But before doing that, we consider the special case 
of Alexandroff topologies, i.e., topologies  with the property that the intersection of any collection of open sets is open. Alexandroff topologies are typical $T_0$ but  not $T_1$ (the discrete topology is the only $T_1$ Alexandroff topology) and are exactly those topologies that are closed as subsets of $2^X$ \cite{todoruzca,UzVi05}. 

\begin{theorem}\cite{todoruzca}
\label{nwd-alex}
Let $\ideal$ be an ideal over a countable set $X$. Then $\ideal=\nwd(\tau)$
for some Alexandroff topology $\tau$ over $X$ if, and only if,
$\ideal$ is  isomorphic to a free sum of ideals belonging to the
following family: principal ideals, $\fin$, $\fin\times\{\emptyset\}$ and
$\nwd({\Q})$.
\end{theorem}

Now we analyze the case when $\tau$ is Hausdorff. It was known that there is no Hausdorff topology $\tau$ such that $\nwd(\tau)= \mbox{FIN}$ (see \cite{CieJas95}). In fact, there is a more general result.

\begin{theorem}\cite{todoruzca} 
\label{compnwd}
Let $\tau$ be an analytic Hausdorff topology over a countable set $X$ without isolated points. Then, 

\begin{itemize}
\item[(i)]  $\nwd(\tau)$ is  $\mathbf{\Pi}^1_2$ and at least  $F_{\sigma\delta}$.

\item[(ii)]  If there is an $\fsig$  set $\base\su 2^X$  which is a base for $\tau$, then $\nwd(\tau)$ is $\pieb$.  
\item[(iii)]  If $(X,\tau)$ is Fr\'echet and regular, then it has a countable $\pi$-base  (see  Shibakov \cite[Corollary 2]{Shiva2017}). Therefore, $\nwd(\tau)$ is   $F_{\sigma\delta}$-complete.
\end{itemize}
\end{theorem}

A typical example of a topology with an $\fsig$ base   is $CL(\cantor)$, the collection of clopen subsets of $\cantor$ with the product topology.  In \cite{Todoruzca2014} it is shown that $\nwd(CL(\cantor))$ is Borel. So, a natural question is 

\begin{question}
Let $(X,\tau)$ be a countable topological space. Suppose $\tau$ has an $\fsig$ base. Is $\nwd(\tau)$ Borel?
\end{question}

In  \cite{MKWELA2018} some nice examples of Hausdorff topologies on $\nat$  are presented, whose nwd ideal has some applications in number theory. 

\begin{example}
\label{nonwd}
\cite{todoruzca}
The following ideals are not of the form $\nwd(\tau)$ for any Hausdorff topology (on the corresponding set).
\begin{itemize}
\item[(i)] $\emptyset\times \ideal$ for any $\fsig$ ideal $\ideal$.

\item[(ii)] The ideal of all subsets of $\omega^2$ with order type smaller than $\omega^2$. 

\item[(iii)] The ideal of scattered subsets of $\Q$ (i.e., subsets of $\Q$ which do not contain an order isomorphic copy of $\Q$).
\end{itemize}
\end{example}

\begin{question}
Find general conditions guaranteeing that a given ideal on a countable set is of the form $\nwd(\tau)$ for $\tau$ a Hausdorff topology. 
\end{question}


\subsection{Topological representation.}
\label{toprepre2}

Suppose $X$ is a Polish space and $J$ is a $\sigma$-ideal of subsets of  $X$. Let $D\su X$ be a countable dense set.  An ideal  $\ideal_J$  on $D$ is defined as follows (see \cite{Sabok2011,AKWELA2015} and the references therein). Let $A\su D$; then, 
\[
A\in \ideal_J\;\; \Leftrightarrow \;\; \cl{A}\in J.
\]
An ideal $\ideal$ on $\nat$ has a {\em topological representation} \cite{AKWELA2015} if there is Polish space $X$, a $\sigma$-ideal $J$ on $X$ and a countable dense set $D\su X$ such that $\ideal$ is isomorphic to $\ideal_J$.  Notice that, by definition, $\nulo(\Q)$ has a topological representation in $\R$. Both ideals $\nwd(\Q)$ and $\nulo (\Q)$ are tall and  $F_{\sigma\delta}$.  In \cite{FarahSolecki2003} it is shown that $\nwd(\Q)$ and $\nulo(\Q)$ are not isomorphic and also that  none  of them is a $P$-ideal.

Topological representable ideals have the following interesting characterization. An ideal $\ideal$ is {\em countably separated} if there is a countable collection 
$\{X_n:\; n\in \nat\}$ such that  for all $A\in \ideal$ and all $B\nin \ideal$, there is $n$ such that $A\cap X_n=\emptyset$ and  $B\cap X_n \nin \ideal$.  This notion was motivated by the results in \cite{Todor96}.

\begin{theorem}
\label{toprepre3}
\cite[Theorem 1.1]{AKWELA2015}
Let $\ideal$ be an ideal on a countable set.  The following are equivalent: 
\begin{itemize}
\item[(i)]  $\ideal$ has a topological representation.
\item[(ii)] $\ideal$ has  a topological representation on $\cantor$ with an ideal $J$ generated by a collection of closed nowhere subsets of $\cantor$.
\item[(iii)] $\ideal$ is tall and countably separated.
\end{itemize}
\end{theorem}


\begin{theorem}
\cite[Corollary 1.5]{AKWELA2015}
If a co-analytic ideal has a topological representation, then it is either $\pieb$-complete or  $F_{\sigma\delta}$-complete.
\end{theorem}

Let us see some examples of ideals which are not topologically representable.

\begin{example}  
\begin{itemize}
\item[(i)]  Let $\ideal$ be an $\fsig$ ideal on $\nat$. We have already mentioned that $\emptyset\times\ideal$ is not of the form $\nwd(\tau)$ for any Hausdorff topology without isolated points  (see Example \ref{nonwd}). Suppose now that $\ideal$ is not tall. It is easy to verify that   $\emptyset\times\ideal$ is not tall and hence it is not topologically  representable.

\item[(ii)]  Consider the ideal $\ideal$ of all subsets of $\omega^2$ of order type smaller than $\omega^2$ (see Example \ref{nonwd}). Then $\ideal$ is tall but it is not countably separated.  The same happens with the ideal $\fin\times\fin$ (see \cite{AKWELA2015}).

\end{itemize}
\end{example}

In \cite[Proposition 4.3]{AKwela2017} it was shown that every countably separated ideal is weakly selective (denoted \ws\ in  \S\ref{prelimi} ), so the following is a natural question. 

\begin{question}
\cite{Sabok2011}
Let $\ideal$ be  a tall,  weakly selective $F_{\sigma\delta}$ ideal. Does $\ideal$ have a topological representation?
\end{question}

For the previous question, one could start with a Farah ideal instead of a $F_{\sigma\delta}$ (see \S\ref{fsig-fsigdelta}).

Since $\Q$ has a countable basis, then $\nwd(\Q)$ is countably separated. On the other hand, $\nwd(CL(\cantor))$ is not weakly selective and therefore it is is not countably separated (see  Example 3.9 in \cite{CamargoUzca2018b}).  Thus a natural question is the following. 

\begin{question}
Let $(X,\tau)$ be a countable Hausdorff space without isolated points.  When is $\nwd(\tau)$ countably separated?  When is it weakly selective?
\end{question}

\subsection{Marczewski-Burstin representations.}

Let $\filter$ be a family of non empty subsets of $X$. The {\em Marczewski ideal} associated to $\filter$ is defined as follows (see \cite{MKWELA2018} and references therein): 
\[
S^0(\filter)=\{A\su X:\; \forall V\in \filter, \exists W\in \filter \;(W\su V\wedge W\cap A=\emptyset)\}.
\]
If $\tau$ is a topology on $X$ and  $\filter$ is a base for $\tau$, then $S^0(\filter)$ is $\nwd(\tau)$. If an  ideal  $\ideal$ is equal to $S^0(\filter)$ for some family of non empty subsets of $X$, then it  is said  that $\ideal$ is  {\em Marczewski-Burstin representable by} $\filter$.  When such $\filter$ can be found countable, it is said that $\ideal$ is  {\em Marczewski-Burstin countably  representable}, which is denoted $\mathcal{MBC}$. 

\begin{example}
\cite[Theorem 4.12]{MKWELA2018}
$\nulo(\Q)$ is $\mathcal{MBC}$. 
\end{example}

It is clear that when $\filter$ is an analytic collection of subsets of a countable set $X$, then $S^0(\filter)$ is at most $\mathbf{\Pi}^1_2$. 
Analogously to what happen with  $\nwd(\tau)$ (see Theorem \ref{compnwd}), if $\filter$ is an $\fsig$ family, then $S^0(\filter)$ is $\pieb$. 

\begin{theorem}
\cite[Theorem 4.4]{MKWELA2018}
(i) Let $\ideal $ be an $\mathcal{MBC}$ ideal. Then  $\ideal$ is  $F_{\sigma\delta}$ and countably separated. 

(ii) If $\ideal$ is countably separated, then there is a $\mathcal{MBC}$ ideal $\cj$ such that $\ideal \su \cj$. 
\end{theorem}

There are two natural properties about  $\filter$  which imply that $S^0(\filter)$ is  tall (see \cite[Theorem 3.6]{MKWELA2018}).

\begin{question}
\cite{MKWELA2018}
Let $f:\nat\to \Q$ be a bijection.  Let 
\[
\cj_c=\{A\su\nat:\; \cl{f(A)}\;\mbox{is countable}\}.
\]
Then $\cj_c$ is an ideal. Is it $\mathcal{MBC}$?

\end{question}

If $(X,\tau)$ is a countable topological  space without isolated points and has  a countable $\pi$-base, then $\nwd(\tau)$ is isomorphic to $\Q$ and clearly is  $\mathcal{MBC}$. Thus we have the following.

\begin{question} Let $(X,\tau)$ be a contable topological space  without isolated points such that $\nwd(\tau)$ is $\mathcal{MBC}$. Is $\nwd(\tau)$ isomorphic to $\nwd(\Q)$?
\end{question}

\section{Ordering the collection of ideals}

One of the main tools for the study of combinatorial properties of
ideals are some  orders (in fact, pre-order) defined on the collection of all ideals: Kat\v{e}tov order $\leq_K$, Rudin-Keisler order $\leq_{RK}$ and Tukey order $\leq_T$. 

Let $\ideal$ and $\cj$ be two ideals on $X$ and $Y$ respectively. We say that $\ideal$ is {\em Kat\v{e}tov below} $\cj$, denoted
$\ideal\leq_K\cj$, if there is a function $f: Y\to X$ such that
$f^{-1}[E]\in\cj$ for all $E\in\ideal$. If $f$ is finite-to-one, then we
write $\ideal\leq_{KB}\cj$ and refer to the  (pre)order $\leq_{KB}$ as the
\emph{Kat\v etov--Blass order}. We say that two ideals $\ideal$ and $\cj$ are {\em Kat\v etov equivalent}, denoted  $\ideal\approx_K \cj$,  if $\ideal\leq_K\cj$ and $\cj\leq_K \ideal$.    Let $Z(\mathcal{I})$ and $Z(\mathcal{J})$  be the  corresponding spaces defined in  \S \ref{prelimi}.  If $f:Y\to X$ is a function we will abuse the notation and consider $f:Z(\mathcal{J})\rightarrow Z(\mathcal{I})$  by letting $f(\infty)=\infty$. If $f$ is a Kat\v etov reduction between $\ideal$ and $\cj$, then  $f:Z(\mathcal{J})\rightarrow Z(\mathcal{I})$ is clearly continuous. Conversely, if there is  $f:Z(\mathcal{J})\rightarrow Z(\mathcal{I})$ continuous with $f^{-1}(\infty)=\{\infty\}$, then  $\mathcal{I}\leq_K \mathcal{J}$.

Let $(D, \leq)$ be a directed ordered set, i.e., for each $x,y\in D$, there is $z\in D$ such that $x, y\leq z$.  A set $A\su D$ is {\em bounded} if there is $x\in D$ such that $y\leq x$ for all $y\in A$. The dual notion to bounded set is that of cofinal set. A set $A\su D$ is {\em cofinal}, if for each $x\in D$, there is $y\in A$ such that $x\leq y$.  Let $D$ and $E$ be two directed orders. A function $f:D\to E$ is called {\em Tukey}, if preimages under $f$ of sets bounded in $E$ are bounded in $D$.  We write $D\leq_T E$ if there is a Tukey function from $D$ to $E$ and we say that $D$ is {\em Tukey reducible} to $E$.  

We shall focus only on the Kat\v etov order as it is crucial for stating some   important open questions. We shall follow the works of Hru\v s\'ak \cite{Hrusak2017}  and Meza \cite{meza-tesis}  (see also \cite{HMTU2017}) which are basic references on this topic. We refer the reader to \cite{solecki2015,solectodor2011,SolecTodor2004} for results on Tukey order.  
The  {Rudin-Keisler} order will be defined in \S\ref{frechet} to state some questions.

\begin{theorem}\cite{Hrusak2017}
Let $\ideal$ and $\cj$ be two ideals on $\nat$. Then, 
\begin{itemize}
\item[(i)] $\ideal\approx_K \fin$ if, and only if $\ideal$ is not tall.
\item[(ii)] If $\ideal\su \cj$, then $\ideal\leq_K\cj$. 
\item[(iii)] if $X\in \ideal^+$, then $\ideal\leq_K \ideal\restriction X$.
\end{itemize}
\end{theorem}

For many combinatorial properties there are ideals (usually Borel ones of a low complexity)
which are critical with respect to the given property, that is,  they are maximal or minimal in the Kat\v etov order $\leq_K$ among all ideals satisfying the property. To illustrate this we present  some examples (see \cite{HMTU2017} for many other similar results). 
A {\em countable  splitting family} for an ideal $\ideal$ on $\nat$  is a countable collection   $\mathcal{X}$ of infinite subsets of $\nat$ such that for every $Y\in\ideal^+$, there is $X\in\mathcal{X}$ such that $| X\cap Y| = | Y\setminus X|=\aleph_{0}$.

\begin{theorem}
\cite{HMTU2017} 
Let $\ideal$ be a tall ideal on $\nat$.  Then, 
\begin{itemize}

\item[(i)] $\nat\to(\ideal^+)^2_2$ if, and only if,
$\mathcal{R}\nleq_K\ideal$, where $\mathcal{R}$ is the random graph ideal.

\item[(ii)] $\ideal$ is a \qpoint-ideal if, and only if, $\ideal\restriction X \ngeq_{KB}\mathcal{ED}_{\fin}$
    for every $\ideal$-positive set $X$.
    
 \item[(iii)] $\ideal$ admits a countable splitting family  if, and only if, $\conv \leq_K \ideal$. 
\end{itemize}
\end{theorem}

The following theorems show some global properties of the Kat\v etov order.

\begin{theorem}
\begin{itemize}
\item[(i)] (H. Sakai \cite{sakai2018}) The family of all analytic P-ideals has a largest element with respect to  $\leq_{KB}$, and thus also with respect to $\leq_K$. 

\item[(ii)] (H. Sakai \cite{sakai2018}) There is an analytic P-ideal $\mathcal J$ such that $\mathcal{I} \leq_{KB} \mathcal{J}$ for all $F_\sigma$ ideal $\mathcal{I}$.  

\item[(iii)](M. Hru\v s\'ak and J. Greb\'ik \cite{GrebikHrusak2018}) There is no Borel tall ideal $\leq_K$-minimal among all Borel tall ideals.

\item[(iv)] (Kat\v etov, see \cite{sakai2018}) There is no Borel ideal which is $\leq_K$-maximum among all Borel ideals.
\end{itemize}
\end{theorem}
There is a result  similar to part (ii) proved by Hru\v s\'ak-Meza \cite{hrusakmeza2013}  showing that there is a universal analytic P-ideal.

Next  results show two  very interesting dichotomies.  The ideals  $\mathcal{R}$, $\mathcal{ED}$, $\mathcal{Z}$ and $\mathcal{S}$ were defined in \S\ref{ejem}.
 
\begin{theorem}
(M. Hru\v s\'ak \cite{Hrusak2017}) (Category Dichotomy) Let $\ideal$ be a Borel ideal. Then either $\ideal\leq_K \nwd(\Q)$, or there is an $\ideal$-positive set $X$ such that $\mathcal{ED}\leq_K\ideal\restriction X$. 
\end{theorem}

\begin{theorem}
(M. Hru\v s\'ak \cite{Hrusak2017}) (Measure Dichotomy) Let $\ideal$ be an analytic $P$-ideal. Then either $\ideal\leq_K \mathcal{Z}$, or there is an $\ideal$-positive set $X$ such that $\mathcal{S}\leq_K \ideal\restriction X$. 
\end{theorem}

\begin{question}\label{questionSole} (M. Hru\v s\'ak \cite{Hrusak2017}) Is $\cR\leq_K \cs$?
\end{question}

As we mentioned above there is no maximum among Borel ideals; however, we have the following.

\begin{question}  (H. Sakai  \cite{sakai2018})  Let $1\leq \alpha<\omega_1$. Is  there a Borel ideal $\ideal$ such that $\cj\leq_K\ideal$ for all $\mathbf{\Sigma}^0_\alpha$ ideal $\cj$?
\end{question}

The following is a fundamental problem. 

\begin{question}
\label{extensionFsigma} (M. Hru\v s\'ak \cite{HMTU2017})
If $\ideal$ is a Borel  tall ideal, then either there is an $\ideal$-positive set $X$ such that $\ideal\restriction X\geq_K\conv$, or there is an $F_{\sigma}$-ideal $\cj$ containing $\ideal$. 
\end{question}

See Theorem \ref{fsigextesion2} for a  partial answer to the previous question. 

\begin{question}
\label{extensionFsigmadelta} (M. Hru\v s\'ak \cite{HMTU2017})
Does every Borel ideal $\ideal$ satisfy that either $\ideal\geq_K \fin\times\fin$, or there is an $F_{\sigma\delta}$-ideal $\cj$ such that
$\ideal\su \cj$?
\end{question}


\section{Ramsey and convergence properties}

In this section we discuss some properties of ideals which have been motivated by properties of convergent sequences and series on $\R$  \cite{Filipowetal2012,Filipow2010,Filipowetal2007,Filipowetal2008}: Bolzano-Weierstrass, Riemann's rearrangement Theorem and convergence in functional spaces. Those properties have a natural connection with Ramsey's theorem. 

We have not included the game theoretic version of Ramsey properties which is indeed a very interesting approach. We refer the reader to the work of Laflamme \cite{laflamme96,laflamme2002}. 

To each ideal there is an associated notion of convergence that we describe hereunder. Let $X$ be a topological space and $\ideal$ an ideal on $\nat$. A sequence $(x_n)_n$ in $X$ is $\ideal$-convergent to $x$, if $\{n\in \nat:\; x_n\nin U\}\in \ideal$ for every open set $U$ of $X$ with $x\in U$.  Notice that $\fin$-convergence is the usual notion of convergence of sequences. 

We recall that an ideal $\ideal$ is called Ramsey at $\nat$ when it satisfies $\nat\to(\ideal^+)^2_2$, and it is called Ramsey when 
$\ideal^+\to(\ideal^+)^2_2$ (see \S\ref{prelimi}).  \footnote{The reader familiar with   \cite{Filipowetal2007,Filipowetal2008} should notice that what they called a Ramsey ideal  (resp.  h-Ramsey) we have called Ramsey at $\nat$ (resp. Ramsey).} 

An ideal $\ideal$ has the \emph{Bolzano--Weierstrass property}, denoted  $\bw$, if for any bounded sequence $\{x_n:n\in\nat\}$ of real numbers there is an $\ideal$-positive set $A$ such that $\{x_n:n\in A\}$ is $\ideal$-convergent. An ideal  $\ideal$ has the \emph{finite Bolzano--Weierstrass property}, denoted  $\fbw$,  if for any bounded sequence $\{x_n:n\in\nat\}$ of real numbers there is an $\ideal$-positive set $A$ such that $\{x_n:n\in A\}$ is convergent. An ideal $\ideal$ is $\mon$ (or \emph{monotone}), if for any sequence $\{x_n:n\in\nat\}$
of real (equivalently rational) numbers there is an $\ideal$-positive set $X$ such that
$\{x_n:n\in X\}$ is monotone  (possibly eventually constant).  We say that $\ideal$ is hereditarely  mononote, denoted {h-{\mon}}, if $\ideal \restriction A$ is $\mon$ for all $A\nin \ideal$. 
 Neither $\nwd(\Q)$ nor $\ideal_d$ satisfy $\fbw$ (see \cite{Filipowetal2007}).

\begin{theorem}\cite[Theorem 3.16]{Filipowetal2008}
Let $\ideal$ be an ideal on $\nat$. The following are equivalent:
\begin{itemize}
\item[(i)] $\ideal\restriction A$ is {$\fbw$} for every $A\nin \ideal$. 
\item[(ii)]  For every  collection $\{A_s: \; s\in\binary\}$ such that $A_\emptyset\nin \ideal$, $A_s=A_{s\widehat{\;\;}0}\;\cup\; A_{s\widehat{\;\;}1}$ and $A_{s\widehat{\;\;}0}\cap A_{s\widehat{\;\;}1}=\emptyset$ for all
$s\in\binary$. There is $B\in\ideal^+$ and $\alpha\in\cantor$ such that $B\su^* A_{\alpha\restriction n}$
for all $n$. 
\end{itemize}
\end{theorem}

The property (ii) above was denoted \wpp in \cite{CamargoUzca2018b}, and property (i) was denoted {h-\fbw} in  \cite{Filipowetal2007,Filipowetal2008}.
The following theorem summarizes several known results in the literature (see \cite{CamargoUzca2018b} for a proof and references).

\begin{theorem}
\label{ppoint} The following  holds for ideals on a countable set. 

\begin{itemize}

\item[(i)]  \pp\ implies \wpp.

\item[(ii)] \qq\ and \wpp\ together is equivalent to  Ramsey.

\item[(iii)] Ramsey implies  \ws.

\item[(iv)]  \wpp\ implies \Sq.

\item[(v)] $\Pm$ is equivalent to  saying that  for every partition $(F_n)_n$ of a set $A\in \ideal^+$ with each piece $F_n$  in $\ideal$,  there is $S\in\ideal^+$ such that $S\su A$ and  $S\cap F_n$ is finite for all $n$.

\item[(vi)] \ws\ is equivalent to \Sq\ together with \qq.

\end{itemize}
\end{theorem}

The usual proof that \fin\ is a {\fbw} ideal shows in fact more: Any \ppoint-ideal is {\fbw}. 

\begin{theorem}
\label{extendedFsig}
\cite[Theorem 3.4 and 4.1]{Filipowetal2007}  Every ideal that can be extended to an $\fsig$ ideal satisfies $\fbw$.
\end{theorem}

\begin{theorem}
\cite[Fact 3.1 and Corollary 3.10]{Filipowetal2008} If an ideal  is  Ramsey at $\nat$, then it satisfies {\mon}, and if it is {\mon}, then {\fbw} holds. Moreover, any  {\mon} analytic $P$-ideal is Ramsey at $\nat$. 
\end{theorem}

\begin{example}
$\ideal_{1/n}$ is $\fbw$ but not {\mon} (see the remark after Corollary 3.10 in \cite{Filipowetal2008}). 
\end{example}

{\fbw} is a Ramsey theoretic property as stated in the following theorems.

\begin{theorem}
\cite[Theorem 3.11]{Filipowetal2008} Let $\ideal$ be  a \qpoint-ideal. Then the following are equivalent:
\begin{itemize}
\item[(i)] $\ideal$ is Ramsey at $\nat$. 
\item[(ii)] $\ideal$ is \mon. 
\item[(iii)] $\ideal$ is \fbw. 
\end{itemize}
\end{theorem}

We have also  a local version of the previous result. 

\begin{theorem}
\cite[Theorem 3.16]{Filipowetal2008} Let $\ideal$ be  an ideal. Then the following are equivalent:
\begin{itemize}
\item[(i)] $\ideal$ is Ramsey. 
\item[(ii)] $\ideal\restriction A$ is {\mon} for every $A\in\ideal^+$. 
\item[(iii)] $\ideal\restriction A$ is {\fbw} for every $A\in\ideal^+$ and $\ideal $ is \qpoint.
\end{itemize}
\end{theorem}

Perhaps one of the most intriguing question is the following. 

\begin{question}
\label{tallRamseyBorel} (Hru\v s\'ak \cite{HMTU2017}) Is there a tall Ramsey  Borel (or analytic)  ideal?
\end{question}

A partial answer to Question \ref{extensionFsigma} is the following. 

\begin{theorem}\label{fsigextesion2} \cite[Proposition 6.5]{Barbarski2013}).
Let $\ideal$ be an analytic $P$-ideal. The following are equivalent. 
\begin{itemize}
\item[(i)] $\conv\not\leq_K\ideal$.
\item[(ii)] $\ideal$ is {\fbw}. 
\item[(iii)] $\ideal$ can be extended to an $\fsig$ ideal. 
\end{itemize}
\end{theorem}

We note that the equivalence of (i) and (ii) was proven in \cite{meza-tesis} (see section 5.1 in \cite{HMTU2017}), and that (ii)  is equivalent to (iii) for analytic $P$-ideals was proven in \cite[Theorem 4.2]{Filipowetal2007}. But the result was formally stated in \cite[Proposition 6.5]{Barbarski2013}).  This motivates a reformulation of Question \ref{extensionFsigma} as follows (see also Theorem \ref{extendedFsig}).

\begin{question}
\cite[Problem 6.1]{Filipowetal2012}  Let $\ideal$ be a tall Borel  $\fbw$ ideal. 
Can $\ideal$  be extended to an $\fsig$ ideal?
\end{question}

\medskip

Now we turn our attention to another classical convergence property that  can be reformulated  in terms of ideals. 
A classical theorem of Riemann says that any conditional  convergent series of reals numbers can be rearranged to converge to any given real number or to diverge to $+\infty$ or $-\infty$.  In other words, if $(a_n)_n$  is a conditional convergent series and $r\in \R\cup\{+\infty, -\infty\}$, there is a permutation $\pi:\nat\to\nat$ such that $\sum_n a_{\pi(n)}=r$. In \cite{Filipow2010,klingaNowik2017} is considered a property of ideals motivated by Riemann's theorem. Let us say that an ideal $\ideal$ has the property $\sf R$, if for any conditionally convergent series $\sum_n a_n$ of real numbers and for any $r\in \R\cup\{+\infty, -\infty\}$, there is a permutation $\pi:\nat\to\nat$  such that $\sum_n a_{\pi(n)}=r$ and
\[
\{n\in \nat:\; \pi(n)\neq n\}\in \ideal.
\]
Similarly, $\ideal$ has property $\sf W$, if for any conditionally convergent series of reals $\sum_n a_n$, there exists $A\in\ideal$ such that the restricted series $\sum_{n\in A}a_n$is still conditionally convergent. In \cite{klingaNowik2017} it is studied similar properties but for series of vectors in $\R^2$. 

\begin{theorem}
(Filip\'ow-Szuca \cite{Filipow2010}) Let $\ideal$ be an ideal on $\nat$. Then, 
\begin{itemize}
\item[(i)] If  $\ideal$  has the property $\sf R$, then it is tall. 
\item[(ii)]  No summable ideal has property $\sf R$. 
\item[(iii)] If $\ideal$ is not $\bw$, then it has property $\sf R$. 
\end{itemize}
\end{theorem}

For instance, since $\ideal_d$ is not $\bw$, then it has property $\sf R$. 

\begin{theorem}
(Filip\'ow-Szuca \cite[Theorem 3.3]{Filipow2010}) Let $\ideal$ be an ideal on $\nat$. The following statements are equivalent.
\begin{itemize}
\item[(i)] $\ideal$  has the property $\sf R$. 
\item[(ii)] There is no a summable ideal $\cj$ such that $\ideal \su \cj$. 
\item[(iii)] $\ideal $ has the property $\sf W$. 
\end{itemize}
\end{theorem}

\begin{question}
(Klinga-Nowik, \cite{klingaNowik2017}) Suppose that (i)  $\ideal$ has the $\sf R$ property; (ii) $\sum_n a_n$ is a conditionally convergent series of reals;
(iii) $\sum_n b_n$ is divergent and all $b_n$ are positive reals. Does there exist $W\in\ideal$ such that $\sum_{n\in W} a_n$ is conditionally convergent and  $\sum_{n\in W} b_n = \infty$?
\end{question}

\medskip

Now we will look at some convergence properties on spaces of continuous functions.  We start with the classical Arzel\'a-Ascoli's theorem characterizing compactness on the pointwise topology.

\begin{theorem}\cite[Theorem 3.1]{Filipowetal2012} (Ideal Version of Arzel\'a-Ascoli Theorem). Let $\ideal$  be an ideal on $\nat$. The following conditions are equivalent.
\begin{itemize}
\item[(i)] $\ideal$ is a  $\bw$ ($\fbw$, respectively).
\item[(ii)] For every uniformly bounded and equicontinuous sequence $(f_n)_{n\in \nat}$ of continuous real-valued functions defined on $[0,1]$, there exists $A\in \ideal^+$  such that $(f_n)_{n\in A}$ is uniformly $\ideal$-convergent (uniformly convergent, respectively).
\end{itemize}
\end{theorem}

Now we present an ideal version of  the classical Helly's selection theorem  in the  space of monotone functions on the unit interval. 

\begin{theorem}\cite[Theorem 5.8]{Filipowetal2012}  (Ideal Version of Helly's Theorem). Let $\ideal$ be an ideal on $\nat$. Suppose that $\ideal$ can be extended to an $\fsig$  ideal. Then for every sequence $(f_n)_{n\in \nat}$ of uniformly bounded monotone real-valued functions defined on $\R$ there is $A\in \ideal^+$ such that the subsequence $(f_n)_{n\in A}$  is pointwise convergent.
\end{theorem}

We recall that,  by Theorem \ref{extendedFsig}, any ideal that can be extended to an $\fsig$ ideal satisfies {\fbw}; thus,  we have the following natural question.

\begin{question}\cite[Problem 5.10]{Filipowetal2012} Let $\ideal$  be an ideal on $\nat$. Are the following conditions equivalent?
\begin{itemize}
\item[(i)] $\ideal$ is   an ${\bw}$ ideal ({$\fbw$} ideal, respectively).
\item[(ii)]  For every uniformly bounded monotone real-valued functions  $(f_n)_{n\in \nat}$ defined on $\R$, there is $A\in \ideal^+$ such that the subsequence $(f_n)_{n\in A}$  is pointwise $\ideal$-convergent (pointwise convergent, respectively).
\end{itemize}
\end{question}

A summary  of  implications among  some of the combinatorial properties studied is as follows. We abbreviate countably generated  and  countably separated by $\omega$-gen and $\omega$-sep, respectively. An ideal is Fr\'echet if it is locally non tall (they will be discussed in the next section).

\bigskip

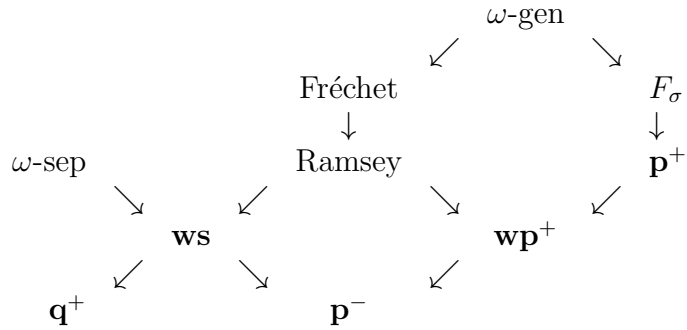
\[
\begin{array}{ccccrcccccccl}
 && &                    &                          & &                  &&  \mbox{} &  &\mbox{$\omega$-gen}& & \\
 && &                    &                          & &                  &  &                         &\swarrow& & \searrow &\\
 && &                    &                          & &             &  &              \mbox{Fr\'echet}              && &  & \fsig\\
 && &                    &                          & &                  &  &        \downarrow                   && &  &\downarrow\\
 && &                    & \mbox{$\omega$-sep}& &                  &  &\mbox{Ramsey}  && &  &\pp\\
&& &    &&\searrow  &&\swarrow& &\searrow&&\swarrow\\
&&&&  &                     &\mbox{\ws}& &&  &  \mbox{\wpp} & \\
&&& & &   \swarrow &&\searrow &&  \swarrow \\
&&&&\mbox{\qq}&& & &\mbox{\Sq}\\
\end{array}
\]


\section{ Fr\'echet ideals}
\label{frechet}

Many of the results presented so far were about tall ideals. In this section we study  Fr\'echet ideals, a  very important class of  non tall  ideals.  This notion  has a topological motivation but  it can be expressed also as a combinatorial  notion.   Recall that to each ideal $\ideal$ on a set $X$ is associate a topological space $Z(\ideal)$ on $X\cup\{\infty\}$  (see \S\ref{prelimi}).  We say that  $\ideal$ is {\em Fr\'echet} if $Z(\ideal)$ is  a Fr\'echet space. Notice  that for  $E\su X$, we have
\[
\mbox{$\infty\in\cl E$ if, and only if,  $E\nin\ideal$.}
\]
It is easy to verify that $\ideal$ is Fr\'echet if, and only if, for every $A\nin\ideal$ there is an infinite $B\su A$ such that every infinite subset of $B$ is not in $\ideal$; that is to say,  $\ideal\restriction A$ is not tall for every $A\nin\ideal$. 
 In other words,  an ideal is Fr\'echet if  it  is locally non tall.

Given a family  $\mathcal A$ of infinite subsets of  $X$,
we define the {\em orthogonal}  of $\ca$ as follows \cite{Todor96}:
\[
\ca^\perp=\{E\su X:\; E\cap A\;\mbox{is finite for all
$A\in\ca$}\}.
\]
Notice that $\ca^\perp$ is an ideal.   If $\ca$ is an analytic family, then  $\ca^\perp$ is co-analytic.

We denote by $\ideal(\ca)$ the ideal generated by $\ca$, that is to say,
\[
E\in \ideal(\ca) \; \Leftrightarrow \; E\su A_1\cup\cdots\cup A_n\;\mbox{for some $A_1, \cdots, A_n \in \ca$}.
\]

\begin{example}
Let $A_n=\{n\}\times\nat$ for $n\in \nat$ and  $\ca=\{A_n:\;n\in\nat\}$. Then 
$\ideal(\ca) = \fin\times\{\emptyset\}$.
\end{example}

The following fact shows the importance of $\perp$ to study Fr\'echet spaces.

\begin{theorem} Let $\ideal$ be an ideal on  $X$. 
\begin{itemize}
\item[(i)]An infinite set  $E\su
X$ is a convergent sequence to $\infty$ in  $Z(\ideal)$  if, and only if, 
$E\in\ideal^\perp$.
\item[(i)] $\ideal$ is Fr\'echet if, and only if,  $\ideal=(\ideal^{\perp})^\perp$.
\end{itemize}
\end{theorem}

\begin{example}
 $(\fin\times\{\emptyset\})^\perp =  \{\emptyset\}\times\fin$ and $(\{\emptyset\}\times\fin)^\perp = \fin\times\{\emptyset\}$. In particular,   $\fin\times\{\emptyset\}$ and $\{\emptyset\}\times \fin$ are Fr\'echet ideals. 
\end{example}

Notice also that $\ca^\perp=((\ca^\perp)^\perp)^\perp$.  In other words, $\ca^\perp$ is a Fr\'echet ideal for any family of sets $\ca$.

 A family  $\ca$ of subsets of  $X$ is {\em almost disjoint} if  $A\cap B$ is finite
for all  $A,B\in\ca$ with  $A\neq B$. Typical examples of  almost disjoint families are the following.

\begin{example}
\label{adf}
(i) For each irrational number $r$, pick a sequence $A_r=\{x^r_n:\;n\in \nat\}$ of rationals numbers converging to $r$. Let $\ca$ be the collection  of all $A_r$ with $r\in \R\setminus \Q$. Then $\ca$ is an almost disjoint family of size $2^{\aleph_0}$.

(ii) Recall that  $2^{<\omega}$ denotes the collection of all finite binary sequences. For each $x\in \cantor$, let $A_x=\{x\restriction n:\;n\in \nat \}$. Then $\{A_x:\; x\in \cantor\}$ is an almost disjoint family. 
\end{example}

As we see next, almost disjoint families are tightly related to Fr\'echet ideals. 

\begin{theorem}\cite{Simon98}
Let $\ideal$ be an ideal on $X$. The following statements are equivalent.

\begin{enumerate}
\item[(i)] $\ideal$ is Fr\'echet.

\item[(ii)] There is an almost disjoint family  $\ca$ of infinite subsets of  $X$ such that
$\ideal =\ca^\perp$.

\item[(iii)]  There is a family  $\ca$ of infinite subsets of  $X$ such that
$\ideal =\ca^\perp$.

\end{enumerate}
\end{theorem}

Let us see  some more examples of Fr\'echet ideals.

\begin{example} 
\label{Ido}
 Consider the ideal $\ideal_{wf}$ generated by the well founded
trees on $\N$ (see \S\ref{complexity}).
The orthogonal of $\ideal_{wf}$ is the ideal $\ideal_{do}$ generated by the
finitely branching trees on $\N$, or equivalently, $\ideal_{do}$
consists of  sets which are dominated by a branch:
\[
A\in \ideal_{do}\;\Leftrightarrow \; \exists \alpha\in \N^\nat \forall
s\in A\forall i<|s| (s(i)\leq \alpha(i)).
\]
The ideal $\ideal_{wf}$ is a complete co-analytic Fr\'echet ideal,  while
the ideal $\ideal_{do}$ is easily seen to be $F_{\sigma\delta}$  (see \cite[Example 2]{dodos2008}). 
\end{example}

\subsection{Selective ideals.}

An ideal is {\em selective}  if it is \ppoint and \qpoint. This is not the original definition given by Mathias  \cite{Mathias77} (who called them  happy families) but it is a reformulation probably due to Kunen.  The original first example of a selective ideal is the following:

\begin{example}
\label{mathias}  (Mathias  \cite{Mathias77}) Let $\ca$ be an analytic almost disjoint family of infinite subsets of $\nat$. Then $\ideal(\ca)$ is a selective ideal.
\end{example}

Next examples were found by Todorcevic \cite{Todor97} in the realm of Banach spaces. 
\begin{example}
\cite[Corollary 7.52]{Todor2010} Let $f, f_n:X\to \R$ be pointwise bounded continuous functions, and suppose that $\{f_n:\; n\in\nat\}$ accumulates to $f$. Let $\ideal_f$ be the ideal defined in Theorem \ref{repre1}, that is to say, 
\[
E\in \ideal_f \;\;\mbox{if, and only if, }  f\nin \overline{\{f_n:\; n\in E\}}.
\]
Then $\ideal_f$ is selective.
\end{example}

One of the reasons for being  interested on selective ideals is due  to the following.

\begin{theorem} (Mathias \cite{Mathias77})
Every selective ideal is Ramsey.
\end{theorem}

Selectivity is the combinatorial counterpart of the topological notion of bisequentiality (see  \cite[Theorem 7.53]{Todor2010})   We only mention the following corollary of this fact which  probably is due to Mathias \cite{Mathias77}.

\begin{theorem}\label{selectivefrechet} Every selective analytic ideal is Fr\'echet.
\end{theorem}

As we already said, if $\ideal$ is analytic, then $\ideal^\perp$ is co-analytic. 
Motivated by the study of Rosenthal compacta  
Krawczyk \cite{Krawczyk92} and Todor\v{c}evi\'{c}  \cite{Todor99,Todor2010} have shown the following (see also  \cite{dodos2008}):

\begin{theorem} If $\ideal$ is a selective analytic ideal
not countably generated, then $\ideal^\perp$ is a complete co-analytic set.
\end{theorem}
 
The  following examples illustrate the previous result. 

\begin{example} 
Let $\ca$ be the almost disjoint family given in Example \ref{adf}(ii), and let $\ideal$ be $\ideal(\ca)$. Then $\ideal$ is selective (see  Example \ref{mathias}) and it is analytic (actually it is $\fsig$), but it is not countable generated.  Hence, $\ideal^\perp$ is $\pieb$-complete (see \cite[Example 1]{dodos2008}). 
\end{example}

\begin{example} 
 In Example \ref{Ido} we presented  the ideal $\ideal_{wf}$ generated by the well founded
trees on $\N$ (see section \ref{complexity}).
The orthogonal of $\ideal_{wf}$ is the ideal $\ideal_{do}$ 
consists of  sets which are dominated by a branch.
The ideal $\ideal_{wf}$ is a complete co-analytic set,  while
the ideal $\ideal_{do}$ is easily seen to be $F_{\sigma\delta}$, it is not countably generated and it is not selective (see \cite[Example 2]{dodos2008}).

\end{example}

\subsection{Orthogonal Borel families.}
\label{OrthogonalBorel}
Two families $\ca$ and $\base$ of infinite subsets of $\nat$ are called {\em orthogonal}, if $A\cap B$ is finite for all $A\in \ca$ and $B\in \base$ \cite{Todor96}.   In this section we are interested in  pairs of orthogonal families which are both Borel.   An example is $\ca=\emptyset\times \fin$ and $\base=\fin\times\emptyset$.  The next theorem says this is the only possible such pair  $(\ideal, \ideal^\perp)$ when one of them is  a $P$-ideal.

 \begin{theorem} (Todor\v{c}evi\'{c}, \cite[Theorem 7]{Todor96}) Let $\ideal$ be an analytic $P$-ideal. Then  $\ideal^\perp$ is countably generated if, and only if,  $\ideal^\perp$ is Borel.
\end{theorem}

 In \cite{GuevaraUzcategui2018}  was constructed a family $\mathfrak{B}$ of $\aleph_1$ non isomorphic Fr\'echet ideals such that both $\ideal$ and $\ideal^\perp$ are Borel.  In fact, every ideal in $\mathfrak{B}$  is $F_{\sigma\delta}$.   Let us recall its definition. 

Let $\{K_n:\; n\in \N\}$ be a partition of $X$.
For $n\in \nat$, let $\ideal_n$ be an ideal on $K_n$. The direct sum,
denoted by $\bigoplus\limits_{n\in \nat} \ideal_n$, is defined  by
\[
A\in \bigoplus_{n\in \N} \ideal_n \Leftrightarrow (\forall n\in \N)(A\cap K_n \in \ideal_n).
\]
For instance, if each $\ideal_n$ is isomorphic to $\fin$, then $\oplus_n\ideal_n$ is isomorphic to $\{\emptyset\}\times \fin$. 
If each $\ideal_n$ is Fr\'echet, then $\bigoplus_n \ideal_n$ is also Fr\'echet. 

The family $\mathfrak{B}$ is the  smallest collection
of ideals on $\N$ containing $\fin$  and
closed under countable direct sums and the operation of taking
orthogonal. The family $\mathfrak{B}$ has some interesting properties.   

\begin{theorem}
\label{teorema_selectivos}
(Guevara-Uzc\'ategui \cite{GuevaraUzcategui2018})
Let $\ideal$ be an analytic selective ideal on $\nat$ and $A\subseteq \N$. The following statements are equivalent:

\begin{itemize}
\item[(i)] $\ideal\restriction A$ is countably generated.
\item[(ii)] $\ideal^\perp\restriction A\in \mathfrak{B}$.
\item[(iii)] $\ideal^\perp\restriction A $ is Borel.
\item[(iv)] $\ideal_{wf}\not\hookrightarrow \ideal^\perp\restriction A$.
\end{itemize}
\end{theorem}

\begin{theorem}
\label{theorem6} (Guevara-Uzc\'ategui \cite{GuevaraUzcategui2018}) For every  $A\su \N^{<\omega}$, the following statements are
equivalent:
\begin{itemize}
\item[(i)] $\ideal_{wf}\restriction A$  belongs to
$\mathfrak{B}$. 

\item[(ii)] $\ideal_{wf}\restriction A$ is
Borel.

\item[(iii)] $\ideal_{wf} \not\hookrightarrow \ideal_{wf}\restriction A$.

\end{itemize}
\end{theorem}

Another interesting co-analytic ideal is  $WO(\Q)$, the collection of  well founded subsets
of $WO(\Q)$. For simplicity, we will write $WO$ instead of
$WO(\Q)$. We first observe that $WO^\perp$ is the ideal of well
founded subsets of $(\Q, <^*)$ where $<^*$ is the reversed order
of $\Q$. In fact, the map $x\mapsto -x$ from $\Q$ onto $\Q$ is an
isomorphism between $WO$ and $WO^\perp$. In particular, $WO$ is a
Fr\'{e}chet ideal.  A linear order $ (L,<)$ is said to be {\em
scattered}, if it does not contain a order-isomorphic copy of
$\Q$.

\begin{theorem}(Guevara-Uzc\'ategui \cite{GuevaraUzcategui2018})
\label{restriction of WO} For every $A\su \Q$, the following statements are
equivalent:
\begin{itemize}
\item[(i)] $A$ is scattered (with the order inherited from $\Q$).

\item[(ii)] $WO\restriction A$ belongs to $\mathfrak{B}$.

\item[(iii)]$WO\restriction A$ is Borel.

\item[(iv)] $WO\not\hookrightarrow WO\restriction  A$.

\end{itemize}
\end{theorem}

\bigskip

It is known that every $\fsig$ tall ideal is not Ramsey,  and also that there is a co-analytic tall Ramsey ideal \cite{HMTU2017}.  We have already stated the basic question of  whether there is a Ramsey tall Borel ideal (see Question \ref{tallRamseyBorel}).  A seemingly weaker question is

\begin{question} Is there a non Fr\'echet Ramsey Borel (or analytic) ideal?
\end{question}

The only Borel Fr\'echet pairs $(\ideal,\ideal^\perp)$ we are aware of are  given by the ideals in $\mathfrak{B}$. So the natural question is:

\begin{question}
\label{ortogonalborel} Is there a Borel Fr\'echet ideal with Borel orthogonal not isomorphic to an ideal in $\mathfrak{B}$?
\end{question}

A related question is the following 
\begin{question}
Are there others  $\pieb$-complete  Fr\'echet ideals satisfying the conclusion of theorem \ref{theorem6}?
\end{question}

Since every Fr\'echet ideal is Kat\v etov equivalent to $\fin$, then Kat\v etov order is trivial among Fr\'echet ideals.  But the  Rudin-Keisler
order is not trivial  on Fr\'echet ideals \cite{GarciaFerreira2013, GarciaOrtiz2017}. 
We say $\ideal\leq_{RK}\cj$ if there is a function $f:\nat\to\nat$ such that
$f^{-1}[E]\in\cj$ if, and only if,  $E\in\ideal$.  

\begin{theorem} (Garc\'ia-Ortiz  \cite{GarciaOrtiz2017})

\begin{enumerate} 
\item[(i)]  There are strictly increasing $\leq_{RK}$-chains of Fr\'echet idelas of size $\mathfrak{c}^+$. Such chains  can be constructed  $\leq_{RK}$-above every Fr\'echet ideal. 

\item[{(ii)}]  For every infinite cardinal $\kappa < \mathfrak{c}$, there is a $\leq_{RK}$-antichain of size $\kappa$.  
\end{enumerate}

\end{theorem}

It is natural then to ask:

\begin{question}
How are the ideals in $\mathfrak{B}$ ordered according to $\leq_{RK}$?
\end{question}

F. Guevara \cite{Guevara2019} has classified the ideals in $\mathfrak{B}$ according to the Tukey order: Except for the countable generated, every ideal in $\mathfrak{B}$ is Tukey equivalent to $\baire$. 

Obviously  a Fr\'echet ideal cannot be topologically representable  as it is  not tall (see Theorem \ref{toprepre3}).  It is easy to check that any Fr\'echet ideal is weakly selective. Thus the following question is appropriate.

\begin{question}
When is a Fr\'echet ideal  countably separated? 
\end{question}

F. Guevara \cite{Guevara2019} has shown that all ideals in $\mathfrak{B}$ are countably separated.

\section{Uniform selection properties}
As we have seen, most of the combinatorial properties for ideals are in fact selection properties.  In this section we analyze the issue of whether the selector can be found Borel measurable. This question can be regarded as one instance of the  classical uniformization problem in descriptive set theory:  Let $B\su X\times Y$ be  a Borel set where $X$ and $Y$ are Polish spaces.  A Borel uniformization for $B$ is a Borel function $F:X\to Y$ such that  $(x,F(x))\in B$ for all $x\in proy_X(B)$. It is well known that, in general, such Borel function does not exist  (see \S 18  of \cite{Kechris94}). 

As an illustration of the problem we are interested, let us consider the notion of tallness.  Let $\cc$ be a  tall Borel (analytic, co-analytic) family of infinite subsets of $\nat$. A very natural question is whether there is a Borel function  $F:\cantor\to \cantor$ such that for all $A\su \nat$ infinite, $F(A)$ is an infinite subset of $A$  and $F(A)\in \cc $.  That is to say, $F$ witness in a Borel way that $\cc$ is tall. In this case we can say that $\cc$ is {\em uniformly tall} or that $\cc$ has a {\em Borel selector}.  This problem was studied in \cite{Grebikuzca2018} and, in particular, they showed that there is a tall $\fsig$ ideal which is not uniformly tall.

\subsection{Uniform Ramsey properties.}
The main question we deal with in this section is whether it is possible to find in a Borel way an homogeneous subset of a given infinite sets.  
This could be briefly stated as whether Ramsey theorem holds uniformly. In the next section we shall see how it can be used to show that a given family is uniformly tall.  Since selective ideals are Ramsey, we start discussing  the uniform versions of the \ppoint and \qpoint\ properties.

We say that  a Borel ideal $\ideal$ is {\em uniformly \ppoint} if there is a Borel function $F$ from $(\cantor)^{\mathbb{N}}$ into $\cantor$ such that whenever $(A_n)_{n}$ is a decreasing sequence of sets in $\ideal^+$, then $A=F((A_n)_{n})$ is in $\ideal^+$ and $A\subseteq^* A_n$ for all $n\in \nat$.
We say that  $\mathcal I$ is {\em uniformly $q^+$}, if there is a Borel function $F$ from $(\cantor)^{\mathbb{N}}$ into $\cantor$ such that whenever $\{K_n\}_n$ is a partition of a set $A$ in $\mathcal{I}^+$ into finite sets, then $S=F((K_n)_n)\subseteq A$, $S$   belongs to $\mathcal{I}^+$ and $|S\cap K_n|\leq 1$ for all $n$.  If $\ideal$ is uniformly \ppoint and \qpoint, we say that $\ideal$ is {\em uniformly selective}.

The following is a uniform version of Theorem \ref{fsigpmas} and Example \ref{mathias}. 

\begin{theorem}\cite{Grebikuzca2018}
\label{unif-p-q} Let $\ideal$ be an $\fsig$ ideal. Then,
\begin{itemize}
\item[(i)] $\ideal$ is uniformly $p^+$.
\item[(ii)] If $\ideal$ is $q^+$, then it is uniformly $q^+$.
\item[(iii)] If $\ca$ is  an almost disjoint family of infinite subsets of $\nat$ which is closed in $\cantor$,  then $\ideal(\ca)$ is uniformly selective.
\item[(iv)] $\fin$ is uniformly selective.
\end{itemize}
\end{theorem}

The previous result naturally suggests the following. 

\begin{question}\cite{Grebikuzca2018} Is $\ideal(\mathcal{A})$ uniformly selective for any  almost disjoint  Borel family $\mathcal A$?  More generally:  is any Borel  selective ideal uniformly selective?
\end{question}

Now we  present some generalization of the Ramsey's theorem.  We need some notation. For $s\in \fin$ and $P\su\nat$ (finite or infinite),  we write $s\sqsubseteq P$ when there is $n\in \mathbb N$ such that
$s=P\cap \{0,1,\cdots, n\}$, and we say that $s$ is an initial segment of $P$.

\begin{theorem}
\label{galvin-lemma} (Galvin's lemma) Let $\filter\su \fin$ and
$M\in\infiN$. There is $N\su M$ infinite such that one of the
following statements holds:
\begin{itemize}
\item[(i)] For all $P\su N$ infinite there is $s\in \filter$ such
that $s\sqsubseteq P$.

\item[(ii)] $N^{[<\infty]}\cap \filter=\emptyset$.
\end{itemize}

\end{theorem}

\noindent Any set $N$ satisfying either (i) or (ii) will be called
 {\em $\filter$-homogeneous}, and the collection of $\filter$-homogeneous sets is denoted by $hom(\filter)$. Notice that if $\filter\su \nat^{[2]}$,
then we have a usual coloring $c:\nat^{[2]}\to \{0,1\}$ by letting $c(s)=1$ if, and only if,  $s\in \filter$. Then, $hom(\filter)=hom (c)$.  Notice also that the previous theorem in particular says that $hom(\filter)$ is a tall family for any $\filter \su\fin$. 

A collection $\mathcal B\su \fin$ is a {\em
front} if it satisfies the following conditions: (i) Every two
elements of $\mathcal B$ are $\sqsubseteq$-incomparable. (ii)
Every infinite subset $N$ of $\nat$ has an initial segment in
$\mathcal B$.  A typical front is $\nat^{[n]}$ for  any $n\in \nat$. 

It is easy to verify that  $hom(\filter)$ is co-analytic subset of $\infiN$ for every
$\filter\su \fin$.  When $\filter\su\base$ and $\base$ is a front,  $hom(\filter)$ is closed in $\infiN$.
We do not know if there is $\filter$ such that $hom(\filter)$ is not Borel.

A key result about the families $hom(\filter)$ is that they are uniformly tall when $\filter\su \base$ for some front $\base$. More precisely:

\begin{theorem}
\cite[Theorem 3.8]{Grebikuzca2018}
\label{uniform-galvin} Let $\base$ be a front. There is a Borel map $S:2^{\base}\times \mathbb{N}^{[\infty]}\rightarrow \mathbb{N}^{[\infty]}$ such that $S(\filter,A)$ is an $\filter$-homogeneous subset of $A$, for all $A\in 
\infiN$  and all $\filter\subseteq \base$.
\end{theorem}

If we use the front $\nat^{[2]}$ we obtain that the classical Ramsey theorem holds uniformly.  We say 
that an ideal $\ideal$ is {\em uniformly Ramsey} if there is a Borel map $S:2^{\nat^{[2]}}\times \mathbb{N}^{[\infty]}\rightarrow \mathbb{N}^{[\infty]}$ such that for all $A\in  \ideal^+$  and all $c: \nat^{[2]}\to \{0,1\}$, $S(c,A)\in \ideal^+$  and it is  a $c$-homogeneous subset of $A$.
The following result is expected.

\begin{theorem}
\cite[Theorem 3.6]{Grebikuzca2018}
Every uniformly selective Borel ideal is uniformly Ramsey.
\end{theorem}

It is also natural to wonder about when  $\fbw$, $\mon$, $\ws$, $\wpp$, etc.  hold uniformly; this is left to the interested reader.

\subsection{Uniformly tall ideals.}
From Theorem \ref{uniform-galvin}, using the front $\nat^{[2]}$,  we obtain that $hom(c)$ is a uniformly tall collection, and  thus $\ideal_{hom(c)}$ is a uniformly tall ideal for any coloring $c$ of pairs of natural numbers.  It should be clear that if a collection $\cc$ contains $hom(c)$ for some coloring $c$, then $\cc$ is also uniformly tall. In fact, most of the examples  we know of uniformly tall families are of that type.   This could be regarded as a method for showing that a given family is uniformly tall (see example \ref{ejemtall} below). 

In particular,  the  random graph ideal $\mathcal R$  (see \S\ref{ejem}) is uniformly tall.  Thus, from the universal property of the random graph, we have  that $\mathcal{R}\leq_K \mathcal{I}$ iff there is a $\filter\subseteq [\mathbb{N}]^2$ such that $hom(\filter)\subseteq \mathcal{I}$. Therefore, if $\mathcal{R}\leq_K \mathcal{I}$, then $\mathcal{I}$ has a Borel selector. That is the case with all examples studied in \cite{Hrusak2017,HMTU2017}.
  Even Solecki's ideal $\mathcal S$ has a Borel selector \cite{GrebikHrusak2018}, even though it is not known whether  it is Kat\v etov above $\mathcal R$ (see Question \ref{questionSole}).  

\begin{example}
\label{ejemtall} 
\cite{Grebikuzca2018} The families of sets listed below are all uniformly tall. This is proved by finding  a coloring $c:\nat^{[2]}\to \{0,1\}$ such that 
$hom(c)$ is a subset of the given family. The coloring used is  the Sierpi\'nski's coloring: Let $X=\{x_n:\;n\in \nat\}$ be a countable set and $\preceq$ a total order on $X$. Define $c:X^{[2]}\to \{0,1\}$ by $c(\{x_n,x_m\})=0$ if, and only if, $n<m$ and $x_n\prec x_m$.   The $c$-homogeneous sets are the $\prec$-monotone sequences in $X$: 
\begin{itemize}
\item[(i)] $\nwd(\tau)$, where  $(X,\tau)$ is a Hausdorff  countable space without isolated points.
\item[(ii)] Let $X$ be a compact metric space and $(x_n)_n$ be a sequence in $X$. Consider 
\[
\cc(x_n)_n= \{A\subseteq \mathbb{N}:\; (x_n)_{n\in A}\;\mbox{is convergent}\}.
\]
\item[(iii)] Let $WO(\mathbb{Q})$ be the collection of all well-ordered subsets of $\mathbb{Q}$ respect the usual order. Let  $WO(\mathbb{Q})^*$  the collection of well ordered subsets of $(\mathbb{Q},<^*)$, where $<^*$ is the reversed order of the usual order of $\mathbb{Q}$. Then, ${\cc}= WO(\mathbb{Q})\cup WO(\mathbb{Q})^*$ is a  tall family. Notice that $\cc$ is $\pieb$-complete. 
\end{itemize}
\end{example}

\bigskip

It is not true  that Galvin's theorem \ref{galvin-lemma} holds uniformly. In fact, there is $\filter \su \fin$ such that $hom(\filter)$ is not uniformly tall (see \cite[Theorem 4.21]{Grebikuzca2018}). Moreover,  there is an $\fsig$ tall ideal which is not uniformly tall (see \cite[Theorem 4.18]{Grebikuzca2018}).  Since the proof of this fact is not constructive, we naturally have the following:

\begin{question} \cite{Grebikuzca2018}
Find a concrete example of an $F_\sigma$ tall ideal without a Borel selector. 
 \end{question}

Tall $\fsig$ ideals are not \qpoint (otherwise they would be selective and thus Fr\'echet, see Theorems \ref{fsigpmas} and \ref{selectivefrechet}). 
This suggests the following:

\begin{question} 
\label{qplusnoselector}
Is there a weakly selective (or \qpoint) tall Borel ideal without a Borel selector?
 \end{question}

Property \qpoint\ might be relevant as the next result suggests. 

\begin{theorem}
Let $\ideal$ be an analytic $P$-ideal. The following assertions are equivalent. 
\begin{itemize}
\item[(i)] $\ideal$ is tall. 
\item[(ii)] $\ideal$ has a continuous selector.
\item[(iii)] $\ideal$ is not \qpoint at $\nat$.
\end{itemize}
\end{theorem}

Since the  generalized summable ideals $\ideal^G_h$ (see \S\ref{summableBanach}) are somewhat similar to $P$-ideals, the previous result naturally suggests  the following. 

\begin{question} 
Let   $\ideal^G_h$ be a generalized summable ideal. Suppose $\ideal^G_h$ is tall. Is it uniformly tall?
 \end{question}

The following result characterizes  tall ideals with  continuous selectors. 

\begin{theorem}(J.~Greb\'ik and M.~Hru\u{s}\'ak \cite[Proposition 25]{GrebikHrusak2018}) Let $\ideal$ be  a Borel tall ideal. Then $\ideal $ has a continuous selector if, and only if, for every  family $\{X_n:\; n\in \nat\}$ of infinite subsets of $\nat$ there is an $A\in \ideal$ such that $A\cap X_n\neq\emptyset$ for all $n\in \nat$.
\end{theorem}

These are the only results concerning the complexity of the selector functions. So we naturally wonder  
if  there is a bound in the Borel complexity of the selector for  Borel tall ideals.

Ideals admitting a topological representation (as defined in \S\ref{toprepre2}) are tall and countably separated. So we have the following question (a negative answer of it will solve Question \ref{qplusnoselector}, as countably separated  ideals are weakly selective \cite[Proposition 4.3]{AKwela2017}).

\begin{question}
Suppose $\ideal$ is a co-analytic ideal with  a topological representation. Is $\ideal$ uniformly tall?
\end{question}

Another question we could ask is whether there is a \textquotedblleft simple basis\textquotedblright\  for the  collection of all tall families. More precisely we have the following question:

\begin{question}
Let $\cc$ be a tall family of infinite subsets of $\nat$. Suppose that $\cc$ is analytic or co-analytic.  Is there $\filter\su \fin$ such that $hom(\filter)\su \cc$?
\end{question}

 The restriction on the complexity is necessary as there is a $\mathbf \Pi^{1}_2$ tall ideal $\ideal$ such that  $hom(\filter)\not\su \ideal $ for all $\filter\subseteq\fin$. In particular, $\ideal$ does not contain any  closed hereditary tall  set (see \cite[Theorem 4.24]{Grebikuzca2018}).

 Some test families  for the previous question are  the following: 
 
 \begin{example}
 \label{testejem}
 (a) Let $\cc_1$ and $\cc_2$ be two tall hereditary families with Borel selector. It is easy to verify that $\cc_1\cap \cc_2$ is also uniformly tall. Let $\base_1$ and $\base_2$ two fronts on $\mathbb{N}$, and $\filter_i\subseteq \base_i$, for $i=0,1$;  is there a front $\base_3$ and $\filter_3\subseteq \base_3$ such that $hom(\filter_3)\subseteq hom(\filter_1)\cap hom(\filter_2)$? Or more generally, given  $\filter_i\su \fin$, for $i=0,1$, is there $\filter_3\su \fin$ such that $hom(\filter_3)\subseteq hom(\filter_1)\cap hom(\filter_2)$? 

 \medskip 
   
(b) Let $\ca$ be an almost disjoint analytic family of infinite subsets of $\nat$. Let $\cc(\ca)$ be $\ideal(\ca)\cup(\ideal(\ca))^\perp$.  Then $\cc(\ca)$ is a $\pieb$ tall family.  The question would be for which families $\ca$ there is $\filter\su\fin$ such that $\hom(\filter)\su \cc(\ca)$.

 \medskip
  
(c) Consider the following generalization of  Example \ref{ejemtall} (ii). Let $K$ be a sequentially compact space, and $(x_n)_n$ be a
sequence on $K$. Let
\[
\cc(x_n)_n= \{A\in\infiN:\; (x_n)_{n\in A}\;\mbox{is convergent}\}.
\]
Then $\cc(x_n)_n$ is tall.

A particular interesting example is for $K$  a separable Rosenthal compacta. By Debs' theorem \cite{Debs1987,Debs2009} (see also
\cite{Dodos2006}), in every Rosenthal compacta, $\cc(x_n)_n$ is
uniformly tall. When $K$ is not first countable $\cc(x_n)_n $ is a complete
co-analytic subset of $\infiN$.  We do not know if there is $\filter\su\fin$ such that $hom(\filter)\su \cc(x_n)_n$. 
\end{example}

\subsection{Uniformly Fr\'echet ideals.}
A Fr\'echet ideal $\ideal$ on a countable set $X$   is {\em uniformly Fr\'echet} if there is a Borel function $f:2^X\to 2^X$ such that for all $A\su X$ with $A\nin \ideal$,  $F(A)\su A$, $F(A)$ is infinite  and $F(A)\in \ideal^\perp$.

\begin{example}
(Guevara \cite{Guevara2011}) All ideals in $\mathfrak{B}$ (see \S\ref{OrthogonalBorel}) are uniformly Fr\'echet. 
\end{example}

In view of the previous result, we have the following variant of Question \ref{ortogonalborel}. 

\begin{question}
(Guevara \cite{Guevara2011}) 
Suppose $\ideal$ is an ideal such that $\ideal$ and $\ideal^\perp$ are Borel and uniformly Fr\'echet. Does $\ideal$ belong to $\mathfrak{B}$?
\end{question}

The definition of a uniformly Fr\'echet ideal does not requiere that it has to be a Borel ideal; however,  we do not have an example of a non Borel uniformly Fr\'echet ideal.  

\begin{example}
(Guevara \cite{Guevara2011}) 
The ideals $\ideal_c$ and $\ideal_{do}$ are both uniformly Fr\'echet Borel ideals and  $\ideal^\perp_c$ and $\ideal^\perp_{do}$  are not uniformly Fr\'echet.
\end{example}

The previous example is a consequence of the following general fact. 
\begin{theorem}
(Guevara \cite{Guevara2011}) 
Let $\ideal$ be a Fr\'echet Borel ideal.  If $\ideal^\perp$ is uniformly Fr\'echet, then $\ideal^\perp$ is Borel.
\end{theorem}

Since Ramsey's theorem holds uniformly (see Theorem \ref{uniform-galvin}), we immediately have the following 

\begin{theorem} Every uniformly Fr\'echet ideal is uniformly Ramsey. 
\end{theorem}
We have seen that every selective analytic ideal is Fr\'echet (see Theorem \ref{selectivefrechet}) and also that  every $\fsig$ selective ideal is uniformly selective. Thus we naturally ask the following:

\begin{question} Is every  uniformly selective $\fsig$ ideal  uniformly Fr\'echet? Or more generally, is every uniformly selective Borel ideal uniformly Fr\'echet?
\end{question}

We have already mentioned in Example \ref{testejem} that $\ideal\cup \ideal^\perp$ is a tall familly for any ideal $\ideal$.  It is easy to check that if $\ideal$ is uniformly Fr\'echet then  $\ideal\cup \ideal^\perp$ is uniformly tall. Thus we have the folllowing.

\begin{question} Let $\ideal$ be a Borel  Fr\'echet ideal such that $\ideal\cup \ideal^\perp$  is uniformly tall. Is $\ideal$ uniformly Fr\'echet?
\end{question}

\textbf{Acknowledgements:}  We would like to thank Francisco  Guevara for the observations he made about the first draft of this paper.

\bibliographystyle{plain}

\end{document}